\documentclass[12pt,reqno,a4paper]{amsart}
\usepackage{
    amsmath,  amsfonts, amssymb,  amsthm,   amscd,
    gensymb,  graphicx, comment,  etoolbox, url,
    booktabs, stackrel, mathtools,enumitem, mathdots,  microtype, lmodern,    mathrsfs, graphicx, tikz,  longtable,tabularx, float, tikz, pst-node, tikz-cd, multirow, tabularx, amscd,  bm, array, makecell, diagbox, booktabs,ragged2e, caption, subcaption }
\usepackage{makecell,slashbox}

\usepackage{xcolor}
\usepackage[utf8]{inputenc}
\usepackage{microtype, fullpage, wrapfig,textcomp,mathrsfs,csquotes,fbb}
\usepackage[colorlinks=true, linkcolor=blue, citecolor=blue, urlcolor=blue, breaklinks=true]{hyperref}
\usepackage[capitalise]{cleveref}
\usepackage{todonotes}
\usetikzlibrary{positioning}
\usetikzlibrary{shapes,arrows.meta,calc}

\usetikzlibrary{arrows}

\newtheorem{theorem}{Theorem}[section]

\newtheorem{corollary}[theorem] {Corollary}
\newtheorem{definition}[theorem]{Definition}
\newtheorem{example}[theorem]{Example}
\newtheorem{lemma}[theorem]{Lemma}

\newtheorem{proposition}[theorem]{Proposition}
\newtheorem{remark}[theorem]{Remark}

\newcommand\R{\mathbb{R}}
\newcommand\Z{\mathbb{Z}}

\newcommand{\la}{\langle}
\newcommand{\ra}{\rangle}

\newcommand{\TC}{\mathrm{TC}}

\newcommand{\ct}{\mathrm{cat}}
\newcommand{\cl}{\mathrm{cl}}

\newcommand{\zl}{\mathrm{zcl}}
\newcommand{\malpha}{\mathrm{M}_{\alpha}}

\newcommand{\mbalpha}{\overline{\mathrm{M}}_{\alpha}}

\newcolumntype{x}[1]{>{\centering\arraybackslash}p{#1}}

\begin{document}
\title{Higher topological complexity of planar polygon spaces having small genetic codes}

\author[S. Datta ]{Sutirtha Datta}
\address{Department of Mathematics, Indian Institute of Science Education and Research Pune, India}
\email{sutirtha2702@gmail.com}
\author[N. Daundkar ]{Navnath Daundkar}
\address{Department of Mathematics, Indian Institute of Science Education and Research Pune, India}
\email{navnath.daundkar@acads.iiserpune.ac.in}
\author[A. Sarkar ]{Abhishek Sarkar}
\address{Department of Mathematics, Indian Institute of Science Education and Research Pune, India}
\email{abhisheksarkar49@gmail.com}


\begin{abstract}
We study the higher (sequential) topological complexity, a numerical homotopy invariant for the planar polygon spaces. For these spaces with a small genetic codes and dimension $m$, Davis showed that their topological complexity is either $2m$ or $2m+1$. We extend these bounds to the setting of higher topological complexity. In particular, when $m$ is power of $2$, we show that the $k$-th higher topological complexity of these spaces is either $km$ or $km+1.$

\end{abstract}
\maketitle
\maketitle

\section{Introduction}\label{sec:intro}
A  \emph{motion planning algorithm} in a path-connected topological space $X$ is defined  as a section of the free path space fibration $\pi:X^I\to X\times X,$ where
$\pi(\gamma)=(\gamma(0),\gamma(1)), $
and  $X^I$ denotes the free path space of $X$, equipped with the compact open topology.
To analyze the complexity of designing a motion planning algorithm for the configuration space $X$ of a mechanical system, Farber \cite{FarberTC} introduced the notion of topological complexity. 
The \emph{topological complexity} of a space $X$, denoted by $\TC(X)$, is defined as the smallest positive integer $r$ for which $X\times X$ can be covered by open sets $\{U_1,\dots, U_r\}$, with each $U_i$ admitting a continuous local section of $\pi$. 
The number $\TC(X)$ represents the minimal number of continuous rules required to implement a motion planning algorithm in the space $X$.
Farber \cite[Theorem 3]{FarberTC} showed that $\TC(X)$ is a numerical homotopy invariant of a space $X$.

In \cite{RUD2010}, Rudyak introduced the higher analogue of topological complexity.
For a path-connected space $X$, consider the fibration $\pi_k: X^I\to X^k$ defined by
\begin{equation}\label{eq: pik fibration}
  \pi_k(\gamma)=\bigg(\gamma(0), \gamma\bigg(\frac{1}{k-1}\bigg),\dots,\gamma\bigg(\frac{k-2}{k-1}\bigg),\gamma(1)\bigg).   
\end{equation}
The \emph{higher topological complexity} of $X$,  denoted by $\TC_k(X)$, is the smallest positive integer $r$ for which $X^k$ can be covered by open sets $\{U_1,\dots, U_r\}$,  such that each $U_i$ admits a continuous local section of $\pi_k$ .
Note that when $k=2$, the higher topological complexity $\TC_k(X)$ coincides with $\TC(X)$.  The above definition makes also sense for $k=1$, but $\TC_1(X)$ is always equal to $1$ (see \cite{RUD2010}). 
Similar to  $\TC(X)$, the invariant  $\TC_k(X)$ is associated with motion planning problems, where the input includes not only an initial and final point but also an additional $k-2$ intermediate points.

An older homotopy invariant of topological spaces, known as the \emph{Lusternik-Schnirelmann (LS) category}, was introduced by Lusternik and Schnirelmann in \cite{LScat}.  
The LS-category of a space $X$, denoted by $\mathrm{cat}(X)$, is the smallest positive integer $r$ such that $X$ can be covered by open subsets $V_1, \dots, V_r$ such that each inclusion $V_i\hookrightarrow{} X$ is null-homotopic.
The following inequalities were established in  \cite{gonzalezhighertc} shows how the LS category and higher topological complexity are related:
\begin{equation}\label{eq: tcn lbub}
\mathrm{cat}(X^{k-1})\leq \TC_k(X)\leq \mathrm{cat}(X^k).
\end{equation}

Determining the precise values of these invariants is usually a challenging task. Over the past two decades, several mathematicians have significantly contributed to approximating these invariants with bounds. 
To be more specific, Farber \cite[Theorem 7]{FarberTC} gave a cohomological lower bound on the topological complexity, and this concept was extended to the higher topological complexity by Rudyak in \cite[Proposition 3.4]{RUD2010}. 
Let $\Delta_k:X\to X^k$ be the diagonal map and  $\zl_k(X)$ denote the cup-length of $\mathrm{ker}(\Delta_k^*)$, where \emph{cup-length} is the length of the largest nontrivial product of cohomology classes. 
Rudyak in \cite[Proposition 3.4]{RUD2010}  proved the following inequality, generalizing Farber's cohomological lower bound on the topological complexity.  
\begin{equation}\label{eq: lb higher tc}
  \TC_k(X)\geq \zl_k(X)+1.   
\end{equation}
We refer to the non-negative integer $\zl_k(X)$ as the \emph{higher zero-divisors-cup-length} of $X$. 
On the other hand, if $\cl(X)$ denotes the cup-length of a cohomology ring of $X$. Then there  is an inequality
\begin{equation}\label{eq: cup-length}
\ct(X) \geq \cl(X)+1  
\end{equation}  
 (see \cite[Proposition 1.5]{CLOT}).
 For a paracompact space, there is a usual \emph{dimensional upper bound} on the higher topological complexity given as follows:
\begin{equation}\label{eq: usual dim ub}
 \TC_k(X)\leq k\cdot\mathrm{dim}(X)+1.   
\end{equation}
Moreover, an additional upper bound for $\TC_k(X)$  is formulated in terms of the homotopy dimension of the space (see \cite[Theorem 3.9]{gonzalezhighertc}).

In this paper, we are interested in studying the (higher) motion planning problem for planar polygon spaces. These spaces can be viewed as equivalence classes of oriented planar $n$-gons with consecutive side lengths $\alpha_1, \dots, \alpha_n \in (0, \infty)$ for some $n \in \mathbb{N}$, identified under translation, rotation, and reflection. Practically, one can regard the sides of a polygon as the linked arms of a robot. Then the higher topological complexity of a polygon space describes the minimum number of motion-planning rules required to maneuver the robot from an initial configuration to a final configuration, ensuring the motion passes through a fixed sequence of intermediate configurations. 

We now briefly recall what planar polygon spaces are and some basic information. A \emph{length vector} is a tuple of positive real numbers. The \emph{planar polygon space} associated with a length vector $\alpha = (\alpha_1, \dots, \alpha_n)$, denoted by $\malpha$, is the collection of all closed piecewise linear paths in the plane with side lengths $\alpha_1, \alpha_2, \dots, \alpha_n$, considered up to orientation-preserving isometries.
Equivalently, we can describe $\malpha$ as 
\[\malpha:= \{(v_{1},v_{2},\dots,v_{n})\in (S^{1})^{n} \mid \sum_{i=1}^{n}\alpha_{i}v_{i} = 0 \}/\mathrm{SO}_{2},\]
where $S^{1}$ is the unit circle and the group of orientation-preserving isometries $\mathrm{SO}_{2}$ acts diagonally.
The planar polygon space (associated with $\alpha$) viewed up to isometries is defined as 
\[\mbalpha\coloneqq \{(v_{1},v_{2},\dots,v_{n})\in (S^{1})^{n} \mid \displaystyle\sum_{i=1}^{n}\alpha_{i}v_{i} = 0 \}/\mathrm{O}_{2},\]
where $\mathrm{O}_2$ is the orthogonal group acting on $\mathbb{R}^2$, including all linear isometries.
  
The spaces $\malpha$ and $\mbalpha$ are also called \emph{moduli spaces of planar polygons}. It is clear that $\malpha$ is a double cover of $\mbalpha$.
A length vector $\alpha$ is called \emph{generic} if $\sum_{i=1}^{n}\pm \alpha_{i} \neq 0$. 
For such a length vector $\alpha$, the moduli spaces $\malpha$ and $\mbalpha$ are closed, smooth manifolds of dimension $m:=n-3$ (see \cite{Farberbook}). 
In the rest of this paper, the length vectors are assumed to be generic unless stated otherwise.  

The planar polygon spaces have been studied extensively. For example, Farber and Sch\"{u}tz \cite{farber2007homology} proved that the integral homology groups of $\malpha$ are torsion-free. They also described the Betti numbers in terms of the combinatorial data associated with the length vector.
The mod-$2$ cohomology ring of $\mbalpha$ was computed by Hausmann and Knutson in \cite{HK1}. 

The program of studying the motion planning problem for planar polygon spaces was initiated by Donald Davis. In his several works \cite{Davis4, davis2017topological, Davis1, Davis2},  he has shown that in most cases $\TC(\mbalpha)$ is either $2m$ or $2m+1$. 
The estimates for $\TC(\mbalpha)$, are primarily determined by the usual dimensional upper bound and cohomological lower bound, respectively, with only a few exceptions. Davis utilized the genetic code description to investigate these bounds, noting that the homeomorphism type of $\mbalpha$ is determined by its genetic code $\alpha$.
In this paper, we initiate the study of the higher topological complexity of planar polygon spaces, focusing on spaces with small genetic codes. We establish higher analogues of Davis's results \cite{davis2017topological} 
through a systematic investigation, leveraging a novel small genetic code description for the mod-$2$ cohomology ring of $\mbalpha$.

\subsection{Structure of the article:} In \Cref{sec:pps}, we begin with the digression of genetic codes and the description of the mod-$2$ cohomology ring of the planar polygon spaces. Then we compute the exact value of the LS category of these spaces in \Cref{prop:cat} and briefly explain the key strategy for our proofs.

In \Cref{sec:size2}, we compute bounds on $\TC_k(\mbalpha)$ when the genetic code of $\alpha$ contains a gene of size $2$. Our results \Cref{thm:monogenic-size2} and \Cref{thm: tck-gencode-having-gene-an} address the cases where $\alpha$ has a monogenic code of size $2$ and where the genetic code includes a gene of size $2$, respectively.
In particular, in \Cref{thm:monogenic-size2} we identify cases for which $\TC_k(\malpha)$ is either $km$ or $km+1$.
In \Cref{R^m non-zero}, we explicitly consider the case of a genetic code having a gene of size $2$ and under some minor conditions we show that $\TC_k(\mbalpha)$ is either $km$ or $km+1$.

\Cref{sec:size3} consists of the cases of genetic code containing a gene of size $3$. In \Cref{Thm-sec4}, we establish that $\TC_k(\mbalpha)$ is either $km$ or $km+1$ under mild assumptions. In general, we provide a novel bound in \Cref{Thm-sec4'}, which extends \cite[Theorem 2.4]{davis2017topological}.

In \Cref{sec: Type2}, we deal with genetic codes of Type $2$. Specifically, our result \Cref{thm:type2} extends \cite[Theorem 3.1]{davis2017topological} to the sequential case.

In \Cref{sec:size4}, we obtain the bounds for $\TC_k(\mbalpha)$ where $\alpha$ having monogenic codes of size $4$. We begin this section by classifying (see \Cref{lem:phi0-odd}) such genetic codes for which the $R^m\neq 0$ (see \Cref{thm: cring} for the cohomology class $R$). This is crucial in identifying the genetic codes of $\alpha$ for which we can have $\TC_k(\mbalpha)$ to be either $km$ or $km+1$. We have achieved this in \Cref{thm:size4-sharpbds}. We end this section by proving our general result  \Cref{thm:size4general}, which is a higher analogue of Davis's results \cite[Proposition 4.3 and Proposition 4.4]{davis2017topological}.

Finally, in \cref{size3-type1}, we take care of two sorts of genetic codes. \Cref{thm-sec7} and \Cref{thm-sec7'}
deal with the higher topological complexity of the genetic codes consisting of two genes, each of size $3$, with the former one being a strong bound and the latter one generalizing the result from section $5$ of \cite{davis2017topological}. Our final result \Cref{thm-sec''} involves genetic codes of $\emph{Type 1}$, a generalization of \cite[Proposition 6.2]{davis2017topological}.

\section{Planar polygon spaces}\label{sec:pps}
In this section, we provide the necessary background on planar polygon spaces, which is essential for stating and proving our upcoming results.

We denote the set $ \{1,\dots, r\}$ by $[r]$, for $r \in \mathbb{N}$. 
There are two important combinatorial objects associated with the length vector $\alpha=(\alpha_{1},\alpha_{2},\dots,\alpha_{n})$. 
\begin{definition}
A subset $I\subset [n]$ is called $\alpha$-short  if 
\[\sum_{i\in I} \alpha_i  < \sum_{j \not \in I} \alpha_j\]  and $\alpha$-long otherwise.
\end{definition}
We may simply use "short" as a shorthand for $\alpha$-short when the context is clear. The collection of short subsets can be quite large. However, there is another combinatorial object associated with length vectors that further compacts the short subset data.
It is important to note that the diffeomorphism type of a planar polygon space does not depend on the ordering of the side lengths. Therefore, we assume that the length vector satisfies $\alpha_{1}\leq \alpha_{2} \leq \dots\leq \alpha_{n}$. 

For a length vector $\alpha$, consider the collection of subsets of $[n]$ :
\[ S_{n}(\alpha) := \{J\subset [n] : \text{ $n\in J$ and $J$ is short}\} \] along with the  partial order $\leq$ defined as follows: 
\[I\leq J~ \text{if}~ I=\{i_{1},\dots,i_{t}\} ~\text{and}~ \{j_{i},\dots,j_{t}\}\subseteq J ~\text{with}~ i_{s}\leq j_{s} ~\text{for}~ 1\leq s\leq t.\]
\begin{definition}\label{gc}
The genetic code of $\alpha$ is the set of maximal elements of $S_{n}(\alpha)$ with respect to the partial order defined above. 
\end{definition}
If $A_{1}, A_{2},\dots, A_{k}$  are the maximal elements of $S_{n}(\alpha)$ with respect to $\leq$ then the genetic code of $\alpha$ is denoted by $\la A_{1},\dots,A_{k}\ra$. Each $A_i$'s is called a \emph{gene} and a subset $A_i\setminus \{n\}$ is called \emph{gee} for each $1\leq i\leq k$ (see \cite{davis2017topological} for more details).

\begin{example}
Let $\alpha=(1,\dots,1,n-2)$ ($n$-tuple) be a length vector. Then the genetic code of $\alpha$ is $\la \{n\}\ra$. Moreover, $\malpha\cong S^{n-3}$ and $\mbalpha\cong \R P^{n-3}$.
\end{example}

Let $G$ be the genetic code of $\alpha$. Then it is easy to see that $G$ uniquely determines the collection $S_{n}(\alpha)$. Consequently, it follows from \cite[Lemma 1.2]{geohausmann} that if the genetic codes of length vectors $\alpha$ and $\beta$ are the same, then the corresponding planar polygon spaces are diffeomorphic.

We now recall the mod-$2$ cohomology algebra $H^*(\mbalpha; \mathbb{Z}_2)$ of $\mbalpha$ in terms of the genetic code information. This result was originally established by Hausmann and Knutson in \cite[Corollary 9.2]{HK1}, and later reinterpreted by Davis (see \cite[Theorem 2.1]{davis2017topological}) as follows.
\begin{theorem}\label{thm: cring}
    The algebra $H^*(\mbalpha; \mathbb{Z}_2)$ is generated by classes $R, V_1, \dots, V_{n-1} \in H^1(\mbalpha; \mathbb{Z}_2)$ subject to the following relations:
\begin{enumerate}
    \item $RV_i+V_i^2$, for $i \in [n-1]$.
    \item $V_S:=\prod\limits_{i \in S} V_i$, unless $S$ is a subgee.
    \item For every subgee $S$ with $|S| \geq n-d-2$,
    \begin{equation} \label{cohomology relation 3}
       \sum\limits_{T \cap S= \emptyset} R^{d-|T|} V_T, 
    \end{equation}
     where $T$ is a subgee.
\end{enumerate}    
\end{theorem}
The following result is an immediate application of the previous theorem, which is crucial in proving some of the results in the next section.
\begin{corollary}[ {\cite{davis2017topological}\label{basis for Hd}}]
For the genetic code $\langle \{a,n \} \rangle$  of $\alpha$, the following set
$$\{R^d, R^{d-1}V_1, \dots, R^{d-1}V_a\}$$
forms a basis for $H^d(\mbalpha; \mathbb{Z}_2)$, for $1 \leq d \leq n-4$.
\end{corollary}

We now set up some notations that will be used throughout the paper. 
Let $p_j:  X^k \rightarrow X$ be the projection onto the $k$-th factor. There is an induced map on the cohomology ring
$$p_j^*: H^*(X; \mathbb{K}) \rightarrow H^*(X; \mathbb{K})^{\otimes k},$$
where $\mathbb{K}$ is a field.
Now considering $X=\mbalpha$ and $\mathbb{K}=\Z_2$ we set the following notions
\begin{enumerate}
    \item $R_j:= p_j^{*}(R)$
    \item $ w_j:= p_j^{*}(V_1)$
    \item $\bar{R}_j:= R_j+ R_{j-1}$
    \item $\bar{w}_j:=w_j+w_{j-1}$
\end{enumerate}
Observe that $\bar{R}_j\in \mathrm{ker}(\Delta_k^{*})$ and $\bar{w}_j\in \mathrm{ker}(\Delta_k^{*})$, where $\Delta_k:X\to X^k$ is the diagonal map.

\begin{remark}\label{k-th coho}\
\begin{enumerate}
\item From Corollary \ref{basis for Hd}, it follows that $\{R^d, R^{d-1}V_1\}$ is linearly independent. Thus,
$$\{R_j^d, R_j^{d-1}w_j\}$$
forms a linearly independent set in the cohomology algebra $H^d(\mbalpha; \mathbb{Z}_2)^{\otimes k}$.   
\item As part-(1) extends \Cref{basis for Hd}, we repeatedly utilize similar canonical approaches to explore the relations among cohomology classes in $H^d(\mbalpha; \mathbb{Z}_2)^{\otimes k}$ for different kinds of small genetic codes, using these as key tools to extend the results of Davis.
\end{enumerate}
\end{remark}

We now compute the LS category of planar polygon spaces and obtain bounds on the higher topological complexity of $\mbalpha$ whose dimension is $m=n-3$. 

\begin{proposition}\label{prop:cat}
Let $\alpha$ be a length vector. Then
$\ct(\mbalpha)=m+1$. Moreover, we have 
\[(k-1)m+1\leq \TC_k(\mbalpha)\leq km+1.\]
\end{proposition}
\begin{proof}
Observe that for a  gee $J$ of maximal size, we have $R^{m-|J|}\prod_{i\in J}V_i$ generates $H^{m}(\mbalpha;\Z_2)$. Therefore, for some $J\subseteq [n-1]$ we have $R^{m-|J|}\prod_{i\in J}V_i\neq 0$ This gives us $\cl(\mbalpha)=m$. Then using \cite[Proposition 1.7]{CLOT} we have $\ct(\mbalpha)=m+1$. 
Now it is easy to see that the product $\prod_{j=1}^r(R_j^{m-|J|}\prod_{i\in J}w_i)$ is nonzero, giving us $\cl(\mbalpha^r)=rm$ for any postive integer $r$. Consequently, $\ct(\mbalpha^r)=rm+1$. The desired inequality then follows from \eqref{eq: tcn lbub}.
\end{proof}

In most cases, to establish lower bounds for $\TC(\mbalpha)$, Davis utilized two key group homomorphisms, $\phi$ and $\psi$. Specifically,
 $\phi: H^m(\mbalpha; \mathbb{Z}_2) \rightarrow \mathbb{Z}_2$ is the Poincar\'e duality isomorphism, and by choosing an uniform homomorphism $\psi:  H^{m-1}(\mbalpha; \mathbb{Z}_2) \rightarrow \mathbb{Z}_2$ such that $(\phi \otimes \psi)(z) \neq 0$, where $z \in H^{2m-1}(\mbalpha \times \mbalpha; \mathbb{Z}_2)$ a product of $(2m-1)$-many cohomology classes from the $\mathrm{ker}(\Delta^*)$.
 Thus, one can conclude that $\zl(X) \geq 2m-1$.
 
For example, if the genetic code of a length vector $\alpha$ is $\la \{a, a+b, n\} \ra$, then  the \emph{uniform homomorphism} in the sense of Davis is a homomorphism $\psi: H^{m-1}(\mbalpha)\to \Z_2$ satisfying 
\begin{itemize}
    \item $\psi(V_i)=\psi(V_j)$ if $i,j\leq a$ or if $a<i,j\leq a+b$, and 
    \item $\psi(V_iV_j)=\psi(V_iV_k)$ if $j,k\leq a$ or if $a<j,k\leq a+b$.
\end{itemize}
(Above, we have omitted writing powers of $R$ accompanying $V$'s as done by Davis in \cite{davis2017topological}).

Our main strategy is to determine lower bounds for $\TC_k(\mbalpha)$ as follows. We will construct  classes $\xi \in H^{2m-1}(\mbalpha \times \mbalpha; \mathbb{Z}_2)^{\otimes l}$ when $k=2l$ and $\xi' \in H^{2m-1}(\mbalpha \times \mbalpha; \mathbb{Z}_2)^{\otimes l} \otimes H^m(\mbalpha; \mathbb{Z}_2)$ when $k=2l+1$ such that
 $(\phi \otimes \psi)^{\otimes l}(\xi) \neq 0$  and $((\phi \otimes \psi)^{\otimes l} \otimes \phi)(\xi') \neq 0$.
Here the map $(\phi \otimes \psi)^{\otimes l}$ is given by the usual tensor product of maps as
$$(\phi \otimes \psi)^{\otimes l}:(H^m(\mbalpha; \mathbb{Z}_2)  \otimes H^{m-1}(\mbalpha; \mathbb{Z}_2))^{\otimes l} \rightarrow \mathbb{Z}_2^{\otimes l} \rightarrow \mathbb{Z}_2,$$
defined by $(\phi \otimes \psi)^{\otimes l}(w_1 \otimes \cdots \otimes w_l)= (\phi \otimes \psi)(w_1) \cdots (\phi \otimes \psi)(w_l)$ for $w_1, \dots, w_l \in H^m(\mbalpha; \mathbb{Z}_2)  \otimes H^{m-1}(\mbalpha; \mathbb{Z}_2)$, and $(\phi \otimes \psi)(w'\otimes w'')= \phi(w') \psi(w'')$ for $w' \in H^m(\mbalpha; \mathbb{Z}_2)$, $w'' \in H^{m-1}(\mbalpha; \mathbb{Z}_2)$. Also, we will use the following homomorphisms
\begin{equation} \label{proj of phi-psi}
  (\phi \otimes \psi)_j^{\otimes l}:= 1 \otimes \cdots \otimes 1 \otimes \phi \otimes \psi \otimes 1 \otimes \cdots \otimes 1,   
 \end{equation}
    where the $2j-1$ th and $2j$ th entries are the maps $\phi$ and $\psi$ respectively, for $1 \leq j \leq l$.

In \cite{davis2017topological}, Davis has explicitly shown that $\TC(\mbalpha) \in \{2m, 2m+1\}$, for the genetic code $\alpha$ having gene sizes from the table below. For genes of size up to $4$, the following table depicts the number of possible genetic codes for the first few values of $n$.

\begin{table}[H]
    \centering
    \begin{tabular}{|c|c|c|c|c|}
    \hline
        Gene Size & $n=5$ & $n=6$& $n=7$& $n=8$ \\
        \hline
        2 & 4 &5&6&7\\
        3&&5&15&21\\ 
        4&&&4&21\\
        3,3&&&15&35\\
        4,3 $\text{Type 1}$&&&8&20\\
        4,3 $\text{Type 2}$&&&\textbf{10}&10\\
        3,3,3 $\text{Type 2}$&&&\textbf{1}&1\\
            4,3,3 $\text{Type 2}$&&&\textbf{14}&14\\
        4,3,3,3 $\text{Type 2}$&&&\textbf{2}&2\\
        anything, 2 &&8&55&559\\
        \hline
    \end{tabular}
    \caption{Number of occurrences.}
    \label{tab:my_label}
\end{table}

A genetic code $\langle S, S'\rangle$ is said to be of $\emph{Type 1}$ if $1 \in S \cap S'$. For the case of $n=7$ and onward, we need the notion of \emph{Type 2}. These are the genetic codes of specified size and which are not of $\emph{Type 1}$ and their gees are those of a genetic code with $n=7$. Note that there are $27$ genetic codes of Type $2$ as highlighted in \Cref{tab:my_label}. We refer the reader to \cite[Table 2]{davis2017topological} for the list of all possible gees of the genetic codes of Type $2$ cases.   For example, $\la \{1,2,4,n\}, \{3,4,n\}, \{2,5,n\}, \{1,6,n\} \ra $ is the genetic code of Type $2$ having genes of size $4$, $3, 3, 3$  (see \cite{davis2017topological} for more details).

\section{The case of genetic code having gene of size 2}\label{sec:size2}
In this section, we obtain sharp lower bounds on  $\TC_k(\mbalpha)$ when the genetic code of $\alpha$ has a gene of size $2$.

Since $\mbalpha$ is a connected smooth manifold of dimension $m=n-3$, we know that  $\TC_k(\mbalpha)\leq km+1$.
In \cite[Theorem 2.2]{davis2017topological}, Davis has shown that the $\TC(\mbalpha)\in\{2m, 2m+1\}$ if the genetic code of $\alpha$ is $\langle \{a,n \} \rangle$. We will now generalize this result in the higher setting. 
\begin{theorem}\label{thm:monogenic-size2}
Let $\langle \{a,n \} \rangle$ be the genetic code of $\alpha$. Then
\begin{equation}\label{eq:tck-gencode-an}
 km-\bigg\lfloor\frac{k}{2}\bigg\rfloor +1\leq  \TC_k(\mbalpha)\leq km+1.   
\end{equation}
Moreover, if $m=2^t$ for some non-negative integer $t$, then 
\begin{equation}\label{eq:tck-gc-an-sharp-lbd}
km\leq\TC_k(\mbalpha)\leq km+1.    
\end{equation}
\end{theorem}
\begin{proof}
For $k = 2$, to obtain the lower bound on $\TC(\mbalpha)$, Davis proved that $\bar{w}_2^m \bar{R}_2^{m-1} \neq 0$ by expanding it in bidegree $(m - 1, m)$ as follows:
\[
\bar{w}_2^m \bar{R}_2^{m-1} = \beta R_1^{m - 2} w_1 R_2^{m - 1} w_2 + R_1^{m - 1} R_2^{m - 1} w_2,
\]
where $\beta = \binom{2m - 1}{m - 1} + 1$, and using the linear independence of $\{R_1^{m - 1}, R_1^{m - 2} w_1\}$ together with the fact that $R_2^{m - 1} w_2 \neq 0$.

We now build on this idea to argue in the general case $k\geq 2$.
We first prove the left inequality of \eqref{eq:tck-gencode-an} when $k=2l$.
Consider the following product of cohomology classes:
$$\bar{A}_{l}:= \prod\limits_{j=1}^{l}(\bar{w}_{2j})^m (\bar{R}_{2j})^{m-1}~.$$
  We now expand $(\bar{w}_{2j})^m (\bar{R}_{2j})^{m-1}$ in the multidegree $(0,\dots,0,m-1,m,0,\dots,0)$, where $m-1$ and $m$ appear at $2j-1$-th and $2j$-th position, respectively. We then have
  \begin{equation}\label{eq: 2jth-product}
 (\bar{w}_{2j})^m (\bar{R}_{2j})^{m-1}=\beta R_{2j-1}^{m-2}w_{2j-1}R_{2j}^{m-1}w_{2j}+ R_{2j-1}^{m-1}R_{2j}^{m-1}w_{2j}.    
  \end{equation}
 By comparing the $2j-1$ th position in \eqref{eq: 2jth-product} and applying Remark \ref{k-th coho}, the above product is non-zero.
Now, $\bar{A}_l$ has the following expression:
$$\bar{A}_l=\prod_{j=1}^l \left(\beta R_{2j-1}^{m-2}w_{2j-1}R_{2j}^{m-1}w_{2j}+ R_{2j-1}^{m-1}R_{2j}^{m-1}w_{2j}\right).$$
Observe that the above expansion contains the following  $(m-1,m,m-1,m,\dots,m-1,m)$-multidegree term $$\prod_{j=1}^{l} R_{2j-1}^{m-1}w_{2j}R_{2j}^{m-1}.$$
Every element in the above product is non-zero, which follows from the \Cref{k-th coho} and the fact $V_1R^{m-1}\neq 0$. Moreover, this expression can not be annihilated by any other term, since no other term in the expansion of $\bar{A}_l$ is devoid of $V_1$ in the odd positions.
 
Thus, 
 $$\zl_{k}(\mbalpha) \geq l(m+m-1)=km-\bigg\lfloor\frac{k}{2}\bigg\rfloor$$ holds, and using the relation (\ref{eq: lb higher tc}) we get the required result.

Next, we obtain the desired result for $k=2l+1$. Consider the following product of cohomology classes:
$$\bar{A}_l \left(\overline{w}_{2l+1}\overline{R}_{2l+1}^{m-1}\right).$$
Observe that the $(m-1,m,m-1,m,\dots,m-1,m)$-multidegree term of the above product is $\bar{A}_lw_{2l+1}R^{m-1}_{2l+1}$, which is nonzero. This gives us the inequality 
$$\zl_{k}(\mbalpha) \geq l(m+m-1)+m=km-\bigg\lfloor\frac{k}{2}\bigg\rfloor$$
giving us the desired lower bound of \eqref{eq:tck-gencode-an}.

Now, we give the sharp lower bound of $\TC_k(\mbalpha)$,
 when $m=2^t$ for some non-negative integer $t$. 
We will now prove the left inequality of \eqref{eq:tck-gc-an-sharp-lbd} for the case $k= 2l$.

 We consider another product of cohomology classes:
$$\bar{C}_{l} = \prod\limits_{j=1}^{l-1}(w_{2j-1}+w_{2j+1}),$$
for $l \geq 2$, and $\bar{C}_1$ is the empty product.
Thus for $l=1$, we have $\bar{A}_1\bar{C}_1=\bar{A}_1\neq0$.

Note that multiplying the class $\bar{A}_{l}$ by the class $(w_{2j-1}+w_{2j+1})$ brings $V_1$ into the $(2j+1)$-th position, keeping it non-zero.
Our claim is
\begin{equation} 
\label{induction lemma} 
 \bar{A}_l \bar{C}_l = (\beta+1) \sum\limits_{j=1 ,  j \text{~odd~}}^{2l} \prod\limits_{i=1, i\neq j}^{2l}w_i \prod\limits_{i=1}^{2l} R_i^{m-1}.   
\end{equation}
We induct on $l$.   Suppose the statement is true for $l=\ell$. Then $\bar{A}_{\ell+1} \bar{C}_{\ell+1}$ is
\begin{align*}
\begin{split}
        &\prod\limits_{j=1}^{\ell+1} \left(\beta  w_{2j}R_{2j}^{m-1}w_{2j-1}R_{2j-1}^{m-2}+w_{2j}R_{2j}^{m-1}R_{2j-1}^{m-1} \right) \prod\limits_{j=1}^{\ell} \left(w_{2j-1} + w_{2j+1}\right)\\
        &= \bar{A}_{\ell}\bar{C}_{\ell} \left(\beta w_{2\ell+2}R_{2\ell+2}^{m-1}w_{2\ell+1}R_{2\ell+1}^{m-2}+w_{2\ell+2}R_{2\ell+2}^{m-1}R_{2\ell+1}^{m-1}\right) \left(w_{2\ell-1}+w_{2\ell+1}\right)\\
        &=(\beta +1) \left(\sum_{j=1, j \text{~odd}}^{2\ell} \prod\limits_{i=1, i\neq j}^{2\ell}w_i\prod\limits_{i=1}^{2\ell}R_i^{m-1}\right)\left(w_{2\ell-1}w_{2\ell+2}R_{2\ell+1}^{m-1}R_{2\ell+2}^{m-1}+w_{2\ell+1}w_{2\ell+2}R_{2\ell+1}^{m-1}R_{2\ell+2}^{m-1}\right)\\
        &=(\beta+1) \sum\limits_{j=1 ,  j \text{~odd~}}^{2\ell+2} \prod\limits_{i=1, i\neq j}^{2\ell+2}w_i \prod\limits_{i=1}^{2\ell+2} R_i^{m-1}.
        \end{split}
\end{align*}
In the penultimate step, we have used the induction hypothesis and the fact that $\beta(\beta+1)$ is even. 
Therefore, whenever $(\beta+1) \neq 0$, it follows that $\bar{A}_l\bar{C}_l \neq 0$, providing us :
\begin{align*}
    \zl_k(\mbalpha) \geq l(m+m-1)+(l-1)=2lm-1=km-1.
\end{align*}
We can impose conditions on $m$ to deduce when $(\beta+1)\neq 0$. For $m=2^t$, Lucas's theorem implies that the binomial coefficient $\beta$ is even. Thus, in this situation, we get 
$$\TC_k(\mbalpha) \geq \zl_k(\mbalpha)+1 \geq km.$$

We now address the case when $k=2l+1$.
Consider the element 
$$D:= \bar{A}_l \bar{C}_l \left(\overline{w}_{2l+1}\overline{R}_{2l+1}^{m-1}\right).$$ 
Using \eqref{induction lemma} we have the expression:
\begin{align*}
    (\beta+1) \sum\limits_{j=1 ,  j \text{~odd~}}^{2l} \prod\limits_{i=1, i\neq j}^{2l}w_i \prod\limits_{i=1}^{2l} R_i^{m-1}  \left(\overline{w}_{2l+1}\overline{R}_{2l+1}^{m-1}\right).
\end{align*}
The rightmost term expands as $(w_{2l+1}+w_{2l})\left( \sum_{r=0}^{m-1} \binom{m-1}{r} R_{2l+1}^rR_{2l}^{m-r-1}\right)$. We consider only the term $w_{2l+1} R_{2l+1}^{m-1}$ from the expansion. Observe that $\bar{A}_l\bar{C}_l\left(w_{2l+1} R_{2l+1}^{m-1}\right)\neq 0$ and therefore $D\neq 0$, as the aforementioned term cannot be annihilated by any other term from the expansion of $D$, due to positional differences. Hence,
\begin{align*}
    \zl_{k}(\mbalpha) \geq l(2m-1)+l-1+1+m-1=m(2l+1)-1=km-1
\end{align*}
Consequently, $\TC_k(\mbalpha) \geq\zl_k(\mbalpha)+1 \geq km,$ when $m= 2^t$ for some $t \in \mathbb{N} \cup \{0\}$.
\end{proof}
\begin{remark} \label{beta}
  For $\beta:= \binom{2m-1}{m-1} +1$, by Lucas's theorem $\beta+1  \neq 0$ $($mod $2)$ when $m=2^t$ for some $t \in \mathbb{N} \cup \{0\}$.   
\end{remark}
For genetic codes having a gene of size two of the length vector $\alpha$, Davis
in \cite[Theorem 2.3]{davis2017topological} has shown that $\TC(\mbalpha)\in \{2m,2m+1\}$. We will now prove a higher analogue of this result. In particular, we have:
\begin{theorem}\label{thm: tck-gencode-having-gene-an}  
Suppose the length vector $\alpha$ has a genetic code with genes $\{a,n\}$. Then 
\[ km-\bigg\lfloor\frac{k}{2}\bigg\rfloor +1\leq \TC_k(\mbalpha)\leq km+1.\]
\end{theorem}
\begin{proof} 
First, we prove our assertion when $k=2l$. Note that, there exists a non-empty subset $\{i_1,\dots,i_t\} \subset [b-1]$ such that $\phi(R^{m-t}V_{i_1} \cdots V_{i_t})=1$, as done in the proof of \cite[Theorem 2.3]{davis2017topological}, where $\phi: H^{m}(\mbalpha;\Z_2)\to \Z_2$ is the Poincar\'e duality isomorphism. Define
$$A_j:=\bar{w}_{2j,i_1}^{m+1-t} \bar{w}_{2j,i_2} \cdots \bar{w}_{2j,i_t} \bar{R}_{2j}^{m-1}, ~\text{where}~ {w}_{j,i_s}:= p_{j}^{*}(V_{i_s}) ~\text{and}~\bar{w}_{j,i_s}={w}_{j,i_s}+{w}_{j-1,i_s} .$$

Now using similar arguments as presented in the proof of \cite[Theorem 2.3]{davis2017topological}, we show that the product $\bar{A}_{l} := \prod_{j=1}^{l} A_j$ is nonzero.
Recall from the proof of \cite[Theorem 2.3]{davis2017topological} that there is a homomorphism $\psi:H^{m-1}(\mbalpha;\Z_2)\to \Z_2$ which sends $R^{m-1}$ and $R^{m-2}V_b$ to $1$ and all other monomials to zero. 

 In the following we apply $(\phi \otimes \psi)^{\otimes l}$ on $\bar{A}_{l}$, using the homomorphisms (\ref{proj of phi-psi}).
\begin{align*}
    \begin{split}
        (\phi \otimes \psi)^{\otimes l}(\bar{A}_l) &= \prod_{j=1}^{l} (\phi \otimes \psi)^{\otimes l}_{j} (A_j) \\
        &= \prod_{j=1}^{l} (\phi \otimes \psi)^{\otimes l}_{j} (\bar{w}_{2j,i_1}^{m+1-t} \bar{w}_{2j,i_2} \cdots \bar{w}_{2j,i_t} \bar{R}_{2j}^{m-1}) \\
        &= \prod_{j=1}^l\phi({w}_{2j-1,i_1}^{m+1-t} {w}_{2j-1,i_2} \cdots {w}_{2j-1,i_t} )\psi({R}_{2j}^{m-1})=1.
    \end{split}
\end{align*}
This implies $\bar{A}_l\neq 0$, thereby completing the proof of our assertion.

Now we consider the case $k=2l+1$. Consider the product: 
$$\bar{A}_l \bar{w}_{2l+1,i_1}\bar{w}_{2l+1,i_2} \dots \bar{w}_{2l+1,i_t}\bar{R}_{2l+1}^{m-t}.$$
Observe that the $(m,m-1,\dots m,m-1,m)$-multidegree term of the above product is $\bar{A}_l{w}_{2l+1,i_1}{w}_{2l+1,i_2} \dots {w}_{2l+1,i_t}{R}_{2l+1}^{m-t}$. Applying $ (\phi \otimes \psi)^{\otimes l} \otimes \phi  $ on the previous product obtains:
\begin{align*}
    \begin{split}
         &\left((\phi \otimes \psi)^{\otimes l} \otimes \phi\right) \left(\bar{A}_l{w}_{2l+1,i_1}{w}_{2l+1,i_2}\dots {w}_{2l+1,i_t}{R}_{2l+1}^{m-t}\right)\\
        &= (\phi \otimes \psi)^{\otimes l}   \left(\bar{A}_l\right) \phi\left({w}_{2l+1,i_1}{w}_{2l+1,i_2}\dots {w}_{2l+1,i_t}{R}_{2l+1}^{m-t}\right)\\
        &=1.
    \end{split}
\end{align*} 
This completes the proof.
\end{proof}

Next we compute sharp bounds on $\TC_k(\mbalpha)$, when the genetic code of $\alpha$ contains a gene $\{a,n\}$.
 In particular, using the result \cite[Example 4.15]{Daundkar-etal-Borsuk-Ulam-PPS} for the genetic code $\langle\{2,4,n\},\{a,n\}\rangle$ we get the following.
\begin{proposition}\label{R^m non-zero}
     Let $a \geq 5$ be odd and $\la\{2,4,n\},\{a,n\}\ra$ be the genetic code of the length vector $\alpha$.
     If $m$ is a $2$-power, then
     \[km\leq\TC_k(\mbalpha)\leq km+1.\]
\end{proposition}
\begin{proof}
 Using \cite[Example 4.15]{Daundkar-etal-Borsuk-Ulam-PPS}   we get that $R^{m}$ is non-zero if $a$ is odd. Thus,
\begin{align} \label{relation on R}
    \begin{split}
        \bar{R}^{2m-1} = &~ \sum_{i=0}^{2m-1} \binom{2m-1}{i} R^i \otimes R^{2m-i-1}\\
        =&~ (\beta + 1) (R^{m-1} \otimes R^m + R^m \otimes R^{m-1}),
    \end{split}
\end{align}
is a non-zero cohomology class (using Remark \ref{beta}), where the class $\bar{R}:=R \otimes 1 + 1 \otimes R \in H^*(\mbalpha \times \mbalpha, \mathbb{Z}_2)$.

Assume $k$ to be even, say $k=2l$. Consider the following product of cohomology classes:
  \begin{itemize}
  \item $\bar{A}_{l}:= \prod\limits_{j=1}^{l}(\bar{R}_{2j})^{2m-1}=(\beta + 1)^l \prod\limits_{j=1}^{l} (R_{2j-1}^{m-1} R_{2j}^m +  R_{2j-1}^{m}  R_{2j}^{m-1})$, for $l \in \mathbb{N}$ (using (\ref{relation on R})).
  \item $\bar{C}_{l} = \prod\limits_{j=1}^{l-1}(R_{2j-1}+R_{2j+1})$ for $l \geq 2$, and $\bar{C}_1$ is the empty product.
\end{itemize}

We inductively show that for each $l \in \mathbb{N}$, the following identity holds.
\begin{equation}\label{power of R}
    \bar{A}_l \bar{C}_l= (\beta +1)^l \left(\sum_{j=1}^{2\ell} (R_j^{m-1} \prod\limits_{i=1, i\neq j}^{2\ell}R_i^m)\right)~.
\end{equation}

Consider the case for $l=1$, the identity $\bar{A}_1 \bar{C}_1=\bar{R}_2^{2m-1}=R_1^{m-1} R_2^m + R_2^{m-1} R_1^m$ holds.

Assuming the required identity holds for some $l = \ell \in \mathbb{N}$, we show that it also holds for the case of $\ell +1$. Note that $\bar{A}_{\ell+1} \bar{C}_{\ell+1}$ expands as
\begin{align*}
\begin{split}
        & \bar{A}_{\ell}\bar{C}_{\ell} \bar{R}_{2 \ell+2}^{2m-1} \left(R_{2\ell-1}+ R_{2\ell+1}\right)\\
        &= (\beta +1)^\ell \left(\sum_{j=1}^{2\ell} (R_j^{m-1} \prod\limits_{i=1, i\neq j}^{2\ell}R_i^m)\right) (\beta + 1) (R_{2\ell+1}^{m-1} R_{2\ell +2}^m +  R_{2\ell+1}^{m}  R_{2\ell+2}^{m-1}) \left(R_{2\ell-1}+ R_{2\ell+1}\right)\\
        &=(\beta +1)^{\ell+1} \left(\sum_{j=1}^{2\ell} (R_j^{m-1} \prod\limits_{i=1, i\neq j}^{2\ell}R_i^m)\right) (R_{2 \ell -1} R_{2 \ell +1}^{m-1} R_{2 \ell +2}^m+ R_{2 \ell -1} R_{2 \ell +1}^{m} R_{2 \ell +2}^{m-1}+ R_{2 \ell +1}^{m} R_{2 \ell +2}^{m})\\
        &=(\beta +1)^{\ell+1} \left(R_{2 \ell+1}^{m-1} \prod\limits_{i=1, i\neq 2 \ell+1}^{2\ell +2}R_i^m+ R_{2 \ell+2}^{m-1} \prod\limits_{i=1, i\neq 2 \ell+2}^{2\ell +2}R_i^m+ \sum_{j=1}^{2\ell} (R_j^{m-1} \prod\limits_{i=1, i\neq j}^{2\ell +2}R_i^m)\right)\\
        &=(\beta +1)^{\ell+1} \left(\sum_{j=1}^{2\ell+2} (R_j^{m-1} \prod\limits_{i=1, i\neq j}^{2\ell+2}R_i^m)\right).
        \end{split}
\end{align*}
Thus, we have proved the identity \eqref{power of R}. Each term in the sum has a different multidegree and contains $R^{m}$ at every position except one, where we have $R^{m-1}$.
This implies $\bar{A}_l\bar{C}_l$ is nonzero whenever $\beta$ is even, more precisely when $m=2^t$ (see \Cref{beta}).
Notice that, to get a lower bound for $\TC_k(\mbalpha)$ where $k= 2l$, we consider the non trivial cohomology class $\bar{A}_l \bar{C}_l$ of length $l(2m-1)+(l-1)=km-1$, and thus  $\TC_k(\mbalpha) \geq km$.

For $k$ odd, say $k=2l+1$, consider the following product $\bar{A}_{l}\bar{C}_{l} \bar{R}_{2\ell+1}^m.$
\begin{align*}
    \begin{split}
        \bar{A}_{l}\bar{C}_{l} \bar{R}_{2l+1}^m = (\beta +1)^l \left(\sum_{j=1}^{2l} (R_j^{m-1} \prod\limits_{i=1, i\neq j}^{2l}R_i^m)\right)\bar{R}_{2l+1}^m.
    \end{split}
\end{align*}
Observe that the product will contain $ (\beta +1)^l \left(\sum_{j=1}^{2l} (R_j^{m-1} \prod\limits_{i=1, i\neq j}^{2l}R_i^m)\right) {R}_{2l+1}^m,$ which is non-zero since its every term is of distinct multidegree and contains $R^m$ at every position except one, where we have $R^{m-1}$.  
\end{proof}

\section{The case of monogenic code of size 3}\label{sec:size3}
In this section, we will derive sharp bounds on $\TC_k(\mbalpha)$ when the genetic code of $\alpha$ is $\langle \{a, a+b, n\} \rangle$.
To proceed, we first recall some relevant notions and results from the proof of \cite[Theorem 2.4]{davis2017topological}.
Building on these, we extend the result to the higher setting.

We now recall the following notations from \cite{davis2017topological}.
\begin{enumerate}
    \item $Y_1:= R^r V_i$ with $i \leq a$,
    \item $Y_2:= R^r V_i$ with $a < i \leq a+b$,
    \item $Y_{1,1}:= R^r V_i V_j$ with $i < j \leq a$,
    \item $Y_{1,2}:= R^r V_i V_j$ with $i \leq a < j \leq a+b$ and
    \item $Y_S$ refers to $R^r V_S$ for an appropriate value of $r$ and $S\subseteq [a+b]$ is a subgee.
\end{enumerate}

Let $\phi_w:= \phi(Y_w)$ and $\psi_{w}:=\psi(Y_w)$ for all possible subscripts $w$ (see \Cref{sec:pps}). Then it was shown in \cite[Proposition 2.5]{davis2017topological} that $\phi$ satisfies the following identities:
$$\phi_{1, 1}= \phi_{1, 2}=1,~ \phi_2= a-1,~ \phi_1= a+ b,~ \phi_0=(a-1)b+ \binom{a-1}{2}.$$

We now briefly recall Davis's strategy to construct the uniform homomorphism $\psi$ and describe two equations that come from the relation \eqref{cohomology relation 3}. 
Note that if $\{i,j\}$ is a subgee of the genetic code $\la\{a,a+b,n\}\ra$, then either $1\leq i<j\leq a$ or $1\leq i\leq a<j\leq a+b$. We denote by $\textbf{R}_{1,1}$, a relation \eqref{cohomology relation 3} corresponding to the subgees of the first kind and by $\textbf{R}_{1,2}$ for the subgees of another kind. 
Then counting subgees (or the cohomology classes of type $Y_1$, $Y_2$, $Y_{1,1}$ and $Y_{1,2}$) of an appropriate type we obtain $\psi(\textbf{R}_{1,1})$ as 

\begin{equation}\label{eq: psi acting on R11}
    \psi_0+ (a-2) \psi_1+ b \psi_2+ \binom{a-2}{2} \psi_{1,1}+ (a-2)b \psi_{1,2}=0
\end{equation}
 and $\psi(\textbf{R}_{1,2})$ as 

 \begin{equation}\label{eq: psi acting on R12}
\psi_0+ (a-1) \psi_1+ (b-1) \psi_2+ \binom{a-1}{2} \psi_{1,1}+ (a-1)(b-1) \psi_{1,2}=0.
 \end{equation}
 
In what follows, we will assume that \emph{$a$ is even.} The other cases will follow similarly. Let $\bar{V}_i=V_i\otimes 1+1\otimes V_i$ for $i=1, a+b$.
Then third equation can be obtained by applying $\phi\otimes \psi$ to the $(m,m-1)$-bidegree expansion of $\bar{V}_1^{2m - 1 - 2^t} \bar{V}_{a+b} \bar{R}^{2^t - 1}$ and equating it to $1$, where $\phi$ is the Poincar\'e duality isomorphism and $2^{t-1}<m\leq 2^t$. The $(m,m-1)$-bidegree expansion of this cohomology class is given by 
\begin{equation*}\label{eq:bidegree expansion when a is even}
 Y_2\otimes Y_1+Y_{1,2}\otimes R^{m-1}+(1+\delta_{m,2^t})R^m\otimes Y_{1,2}+Y_1\otimes Y_2, 
\end{equation*}
where $\delta$ denotes the Kronecker delta.
Then $\phi\otimes \psi$ acting on the above expansion and equating it to $1$ gives the third equation:
\begin{equation}\label{eq: phipsi acting on expansion}
\psi_0+\psi_1+ b\psi_2+\varepsilon\psi_{1,2}=1,
\end{equation}
where $\varepsilon\in\Z_2$ is an irrelevant quantity in solving the system of linear equations formed by \eqref{eq: psi acting on R11}, \eqref{eq: psi acting on R12}, and \eqref{eq: phipsi acting on expansion}  with indeterminate $\psi$.
Davis has shown the existence of $\psi$ by proving that the system of linear equations described above has a solution. This shows the non-zeroness of the cohomology class $\bar{V}_1^{2m - 1 - 2^t} \bar{V}_{a+b} \bar{R}^{2^t - 1}$ and consequently gives us the desired lower bound on $\TC(\mbalpha)$ when the genetic code of $\alpha$ is $\la\{a,a+b,n\}\ra$
with $a$ even.

We now characterize the values of $a$ and $b$ for which the $R^m$ is non-zero, where $m$ is the dimension of $\mbalpha$.
\begin{lemma}\label{a mod lemma} 
Let $\la\{a, a+b,n\}\ra$ be the genetic code of $\alpha$ with $n>a+b>a>0$ and $n\geq 6$. 
Then $R^m \neq 0$ if and only if either of the following cases holds:
\begin{enumerate}
 \item $a \equiv 3 ~(mod~ 4)$,
 \item $a \equiv 0 ~(mod ~4)$ and $b$ even,
 \item $a\equiv 2 ~(mod ~4)$ and $b$ odd.
 \end{enumerate}
\end{lemma}
\begin{proof}
Recall that $\phi(R^m)=\phi_0$, where $\phi$ is the Poincar\'e duality isomorphism. Now
$$\phi_0=   (a-1)b + \frac{(a-1)(a-2)}{2}.$$
    \begin{enumerate}
        \item Suppose $a \equiv 3 ~(mod~ 4)$. Then $(a-1)b$ is even and $ (a-1)(a-2)$ has only one $2$ as a factor. Hence $\frac{(a-1)(a-2)}{2}$ is odd.
        \item Suppose $a \equiv 0 ~(mod ~4)$ and $b$ even. Then $(a-1)b$ is even and $\frac{(a-1)(a-2)}{2}$ is odd.
        \item Since $4$ divides $(a-2)
        $, $\frac{(a-1)(a-2)}{2}$ is even. Also $b$ being odd makes $(a-1)b$ odd.
    \end{enumerate}
    It is easy to check that for all other combinations of values of $a$ and $b$, $\phi_0=0$.
\end{proof}

As a consequence of \Cref{a mod lemma}, we obtain the following sharp bounds.

\begin{theorem}\label{Thm-sec4}
Let $m=2^t$ for some non-negative integer $t$ and $\la\{a, a+b,n\}\ra$ be the genetic code of $\alpha$ such that $a$ and $b$ satisfy either of the conditions of \Cref{a mod lemma}. Then $$km \leq \TC_k(\mbalpha)\leq km+1.$$
\end{theorem}
\begin{proof}
Since $R^m$ is non-zero by the conditions of \Cref{a mod lemma}, the proof follows using similar techniques as used in proving  \Cref{R^m non-zero}.
\end{proof}
Next we use a higher analogue of the choice of the cohomology classes in \cite[Theorem 2.4]{davis2017topological} to prove the following result.

\begin{theorem}\label{Thm-sec4'}
Let $\la\{a, a+b,n\}\ra$ be the genetic code of $\alpha$. Then 
\begin{equation}\label{eq: tck-size3-monogenic}
km-\bigg\lfloor\frac{k}{2}\bigg\rfloor +1\leq \TC_k(\mbalpha)\leq km+1.    
\end{equation}
\end{theorem}
\begin{proof} Suppose that $2^{t-1} < m \leq  2^t$.
Before we proceed, we set some notations that will be used throughout the proof. Define $x_{i}:= p_{i}^{*}(V_{a+b})$, $\bar{x}_{i}=x_i+x_{i-1}$  and  $_iY_w:= p^*_i(Y_w)$ for $1\leq i \leq k$.
We proceed with the following cases.

\noindent\textbf{Case-1}: \emph{Suppose $a$ is even.}

\noindent 
    Assume $k=2l$. Consider the element $$\bar{A}_{l}:= \prod\limits_{j=1}^{l}A_j,$$ 
    where $A_j=\bar{w}_{2j}^{2m-1-2^t}\bar{x}_{2j}\bar{R}_{2j}^{2^t-1}$.
    We apply $(\phi \otimes \psi)^{\otimes l}$ on the $(m,m-1,m, m-1,\dots, m, m-1)$-multidegree term in the expansion of $\bar{A}_{l}$. Recall the notation $(\phi \otimes \psi)^{\otimes l}_i$ from \eqref{proj of phi-psi}.
    Then the identity $(\phi \otimes \psi)^{\otimes l}_i (A_i)=\psi_0+\psi_1+ b\psi_2+\varepsilon \psi_{1,2}$ holds for all $i$ (see \eqref{eq: phipsi acting on expansion}).
   Thus, we obtain the following expression:
$$(\phi \otimes \psi)^{\otimes l} (\bar{A}_{l})= \prod_{j=1}^{l}(\phi\otimes\psi)^{\otimes l}_j (A_j)=(\psi_0+\psi_1+ b\psi_2+\varepsilon\psi_{1,2})^{l}.$$

Since $\eqref{eq: phipsi acting on expansion}$ has a solution as shown by Davis, 
 $(\phi \otimes \psi)^{\otimes l} (\bar{A}_{l})=1$ has also has a solution.
Therefore, $\bar{A}_{l}\neq 0$ for $l \in \mathbb{N}$.
This gives us the left inequality of \eqref{eq: tck-size3-monogenic}.

We use similar idea to show that the product $\bar{A}_{l} \bar{x}_{2l+1}\bar{R}_{2l+1}^{m-1}$ is nonzero when $k=2 l +1$.

Observe that the type of the $(m,m-1,\dots,m,m-1,m)$-multidegree term of  $\bar{A}_{l} \bar{x}_{2l+1}\bar{R}_{2l+1}^{m-1}$ is $\bar{A}_{l}~ _{2l+1}Y_{2}$.
Since $\phi_2=a-1=1$ and $(\phi\otimes \psi)^{\otimes{l}}(\bar{A}_{l})=1$, we get $$((\phi\otimes \psi)^{\otimes{l}}\otimes \phi)(\bar{A}_{l}~ _{2l+1}Y_{2})=(\phi\otimes \psi)^{\otimes{l}}(\bar{A}_{l})\phi( _{2l+1}Y_{2})=1.$$
This proves $\bar{A}_{l} \bar{x}_{2l+1}\bar{R}_{2l+1}^{m-1}\neq 0$ and thus we get the left inequality of \eqref{eq: tck-size3-monogenic}. 

\medskip

\noindent\textbf{Case-2}: \emph{ $a$ is odd and $m \neq 2^{t-1}+1$.} 

\noindent Then for $k=2$, Davis proved that $\bar{V}_1^{m-1}\bar{V}_{a+b}^2\bar{R}^{m-2}\neq 0$ by showing the existence of uniform homomorphism $\psi$, and then applying $\phi\otimes \psi$ on the expansion of $\bar{V}_1^{m-1}\bar{V}_{a+b}^2\bar{R}^{m-2}$ in the bidegree $(m,m-1)$. We generalize this idea for general $k\geq 2$.

Assume $k=2l$. Consider the element $$\bar{A}_{l}:= \prod\limits_{j=1}^{l}A_j,$$ 
    where $A_j=\bar{w}_{2j}^{m-1}\bar{x}_{2j}^2\bar{R}_{2j}^{m-2}$.
    We apply $(\phi \otimes \psi)^{\otimes l}$ on the $(m,m-1,m, m-1,\dots, m, m-1)$- multidegree term in the expansion of $\bar{A}_{l}$. Similar to Case-1, we obtain the following expression:
$$(\phi \otimes \psi)^{\otimes l} (\bar{A}_{l})= \prod_{j=1}^{l}(\phi\otimes\psi)^{\otimes l}_j  (A_j)=(\psi_1+ m(b+1)\psi_2+\varepsilon'\psi_{1,2})^{l}.$$  
In \cite[Theorem 2.4]{davis2017topological}, Davis has shown that $\psi_1+ m(b+1)\psi_2+\varepsilon'\psi_{1,2}=1$ has solution. Here, the value of $\varepsilon'$ is again irrelevant. Consequently, $(\phi \otimes \psi)^{\otimes l} (\bar{A}_{ l})=1$ has a solution.
Therefore, $\bar{A}_{l}\neq 0$ for $l \in \mathbb{N}$.
This gives us the left inequality of \eqref{eq: tck-size3-monogenic}.

We use a similar idea to show that the product $\bar{A}_{l} \bar{w}_{2l+1}\bar{x}_{2l+1}\bar{R}_{2l+1}^{m-2}$ is nonzero when $k=2l +1$. 
Observe that the type of the $(m,m-1,\dots,m,m-1,m)$-multidegree term of  $\bar{A}_{l}\bar{w}_{2l+1}\bar{x}_{2l+1}\bar{R}_{2l+1}^{m-2}$ is $\bar{A}_{l}~ _{2l+1}Y_{1,2}$.
Since $\phi_{1,2}=1$ and $(\phi\otimes \psi)^{\otimes{l}}(\bar{A}_{l})=1$, we get $$((\phi\otimes \psi)^{\otimes{l}}\otimes \phi)(\bar{A}_{l}~ _{2l+1}Y_{1,2})=(\phi\otimes \psi)^{\otimes{l}}(\bar{A}_{l})\phi( _{2l+1}Y_{1,2})=1.$$
This proves $\bar{A}_{l} \bar{w}_{2l+1}\bar{x}_{2l+1}\bar{R}_{2l+1}^{m-2}\neq 0$ and thus we get the left inequality of \eqref{eq: tck-size3-monogenic}. 

\smallskip

\noindent\textbf{Case-3}: \emph{$a$ is odd and $m= 2^{t-1}+1$.} 

\noindent Then for $k=2$, Davis proved that $\bar{V}_1^m\bar{V}_{a+b}^{m-1}\neq 0$ by constructing a uniform homomorphism $\psi$ and applying $\phi\otimes \psi$ to the expansion of $\bar{V}_1^m\bar{V}_{a+b}^{m-1}$ in the bidegree $(m,m-1)$. We use this idea to prove our assertion for $k\geq 2$.

Assume $k=2l$. Consider the element $$\bar{A}_{l}:= \prod\limits_{j=1}^{l}A_j,$$ 
    where $A_j=\bar{w}_{2j}^{m}\bar{x}_{2j}^{m-1}$.
    We apply $(\phi \otimes \psi)^{\otimes l}$ on the $(m,m-1,m, m-1,\dots, m, m-1)$- multidegree term in the expansion of $\bar{A}_{l}$. Similar to the previous cases, we obtain the following expression:
$$(\phi \otimes \psi)^{\otimes l} (\bar{A}_{ l})= \prod_{j=1}^{l}(\phi\otimes\psi)^{\otimes l}_j (A_j)=(\psi_1+ (b+1)\psi_2)^{l}.$$ 

In \cite[Theorem 2.4]{davis2017topological}, Davis has shown that $\psi_1+ (b+1)\psi_2=1$ has solution. Consequently, $(\phi \otimes \psi)^{\otimes l} (\bar{A}_{l})=1$ has a solution.
Therefore, $\bar{A}_{l}\neq 0$ for $l \in \mathbb{N}$.
This gives us the left inequality of \eqref{eq: tck-size3-monogenic}.
 The case when $k$ is odd is exactly the same as the odd case of Case 2.
\end{proof}

\section{The genetic codes of type 2}\label{sec: Type2}
In this section, we study the higher topological complexity of planar polygon spaces having genetic codes of Type $2$. We refer the reader to \Cref{sec:pps} for a brief description of these genetic codes.

It turns out that there are exactly $27$ Type $2$ genetic codes (see \cite[Table 2]{davis2017topological} for more details). To avoid repetition, we have opted not to include \emph{Table 2} from Davis's paper. For these genetic codes, Davis has computed the zero-divisors-cup lengths. We now briefly describe Davis's strategy. 

Suppose $m=n-3\geq 4$ is the dimension of $\mbalpha$. The genetic codes are distributed into two cases:

\noindent   \emph{Suppose $m-1$ is not a  $2$-power or if this is the last case of the \cite[Table 2]{davis2017topological} which is the Type $2$ genetic code $\la \{3,4,n\}, \{2,5,n\}, \{1,6,n\} \ra$}. 
Let $\bar{V}_i=V_i\otimes 1+1\otimes V_i$ for $i=1, 2, 3$. In this case, Davis proved that 
$\bar{V}_1^{m-1} \bar{V}_{2}^2 \bar{V}_3 \bar{R}^{m-3} \neq 0$ by constructing a uniform homomorphism $\psi$, and then applying $\phi \otimes \psi$ to the expansion of the cohomology class in bidegree $(m, m - 1)$, where $\phi$ denotes the Poincar\'e duality isomorphism. The existence of $\psi$ was shown by solving a system of linear equations. 
We now explain his ideas briefly. The homomorphism $\phi\otimes \psi$ acting on the expansion of $\bar{V}_1^{m-1} \bar{V}_{2}^2 \bar{V}_3 \bar{R}^{m-3}$ in bidegree $(m,m-1)$ gives us the following:
\begin{equation}\label{eq: phi psi action type2}
\begin{aligned}
& (\phi_{2,3}+\phi_{1,2,3}) \psi_1 
+ \phi_{1,3}(\psi_2+\psi_{1,2})
+ (m-1)\,\phi_1(\psi_{2,3}+\psi_{1,2,3}) 
\\
& \quad + 
\begin{cases}
   \phi_1 \psi_{1,2,3} & \text{if $m$ is a 2-power},\\[6pt]
   \phi_{1,2,3}\psi_1+ \phi_{1,3}\psi_{1,2} & \text{if $m-1$ is a 2-power},\\[6pt]
   0 & \text{otherwise}.
\end{cases}
\end{aligned}
\end{equation}

Davis first computes the values of $\phi_S$, where $\phi$ is a Poincaré duality isomorphism and $S$ is a subgee of a particular type. 
The key ingredient to achieve this in Davis's method is to construct a matrix whose columns represent all subgees, including $\emptyset$, and rows represent all subgees except $\emptyset$. This matrix is actually a binary matrix in which the entry $1$ appears in the places where the corresponding subgees for the row and column are disjoint.
Each row corresponding to the subgee of this matrix, in fact, represents the relation \ref{cohomology relation 3}. The reader is referred to \cite{27cases, Mapleresults} for an illustration. 
After solving this system of linear equations using \texttt{MAPLE}, Davis obtained the values of $\phi_S$ and observed that $\phi_{1,2,3}=1$, while $\phi_{1,3}=\phi_{2,3}=0$ in the first $26$ cases of the \cite[Table 2]{davis2017topological}.  

A similar approach figures out the values of $\psi_S$, except that the rows correspond only to subgees with more than one element, due to the $|S| \geq n-d-2$ constraint in \eqref{cohomology relation 3}.  
The uniform homomorphism $\psi$ in each case given in \cite[Table 2]{davis2017topological} satisfies $\psi_1=1$ and $\psi_{2,3}=\psi_{1,2,3}=0$. One can see that in this case \eqref{eq: phi psi action type2} equates to $1$  for the first $26$ genetic codes of Type $2$ whose gees are described in \cite[Table 2]{davis2017topological}. For the final case in \cite[Table 2]{davis2017topological}, again using \texttt{MAPLE} 
one can see that there is a uniform homomorphism $\psi$ for which the only nonzero values are $\psi_1$ and $\psi_{1,6}$. Moreover, $\phi_{i,j}=1$ for all subgees $\{i,j\}$ and $\phi_{1,2,3}=0$ as $V_{1,2,3}=0$. Since $\phi_{2,3}\psi_1=1$, again  \eqref{eq: phi psi action type2} equates to $1$. Thereby proving the nonzeroness of $\bar{V}_1^{m-1} \bar{V}_{2}^2 \bar{V}_3 \bar{R}^{m-3}$.

 \emph{ If $m-1$ is a $2$-power and it is not the last case of \cite[Table 2]{davis2017topological}} can be dealt similarly to the previous case. So we don't repeat the explanation of his strategy.

\begin{theorem}\label{thm:type2}
Let $\alpha$ has the genetic code of type $T_2$ in the \Cref{tab:my_label} and $m\geq 4$, then: 
\begin{equation}\label{eq: tck-type2}
       km-\bigg\lfloor \frac{k}{2} \bigg\rfloor +1 \leq \TC_k(\mbalpha) \leq km+1.
   \end{equation}
   \end{theorem}
   \begin{proof}
   Recall from \Cref{sec:pps} that $w_j=p_j^*(V_1)$. Similarly, we define $u_j:=p_j^*(V_2)$ and $v_j:=p_j^*(V_3)$. Similar to Davis, we distribute the genetic codes into two cases.
   
\noindent\textbf{Case 1}: \emph{Suppose $m-1$ is not a $2$-power or $\alpha$ has gees as in the final case in \cite[Table 2]{davis2017topological}.}

Consider the element $$\bar{A}_{l}:= \prod\limits_{j=1}^{l}A_j,$$ 
    where $A_j=\bar{w}_{2j}^{m-1} \bar{u}_{2j}^2\bar{v}_{2j} \bar{R}_{2j}^{m-3}$.
    We apply $(\phi \otimes \psi)^{\otimes l}$ on the $(m,m-1,m, m-1,\dots, m, m-1)$-multidegree term in the expansion of $\bar{A}_{l}$. When $m$ is a power of $2$, we obtain the following:
    \begin{align*}
        \begin{split}
            &(\phi \otimes \psi)^{\otimes l} (\bar{A}_{l})= \prod_{j=1}^{l}(\phi\otimes\psi)^{\otimes l}_j (A_j)\\
            &=((\phi_{2,3}+\phi_{1,2,3})\psi_1+\phi_{1,3}(\psi_2+\psi_{1,2})+(m-1)\phi_1(\psi_{2,3}+\psi_{1,2,3})+\phi_1\psi_{1,2,3})^{l}. 
        \end{split}
    \end{align*} 
    Then the above expression is equal to $1$, since $\phi_{1,2,3}=\psi_1=1$ as explained at the beginning of this section.
    
    Next, when $m-1$ is a power of $2$ (and the gee is of the final type in \cite[Table 2]{davis2017topological}), we obtain the following:
     \begin{align*}
        \begin{split}
            &(\phi \otimes \psi)^{\otimes l} (\bar{A}_{l})= \prod_{j=1}^{l}(\phi\otimes\psi)^{\otimes l}_j (A_j)\\
            &=((\phi_{2,3}+\phi_{1,2,3})\psi_1+\phi_{1,3}(\psi_2+\psi_{1,2})+(m-1)\phi_1(\psi_{2,3}+\psi_{1,2,3})+\phi_{1,2,3}\psi_1+\phi_{1,3}\psi_{1,2})^{l}. 
        \end{split}
    \end{align*}
     For the $27$-th gee in \cite[Table 2]{davis2017topological}, both the above two expressions ($m-1$ is a $2$-power and not a $2$-power) become $1$, since there is homomorphism $\psi$ whose only non-zero values are $\psi_1$ and $\psi_{1,6}$. Moreover, $\phi_{i,j}=1$ for all $\{i,j\}$ and $\phi_{1,2,3}=0$. Therefore, only $\phi_{2,3}\psi_1$ survives. 
    
    Finally, in all other situations under this case, we get the following:
    \begin{align*}
        \begin{split}
            &(\phi \otimes \psi)^{\otimes l} (\bar{A}_{l})= \prod_{j=1}^{l}(\phi\otimes\psi)^{\otimes l}_j (A_j)\\
            &=((\phi_{2,3}+\phi_{1,2,3})\psi_1+\phi_{1,3}(\psi_2+\psi_{1,2})+(m-1)\phi_1(\psi_{2,3}+\psi_{1,2,3}))^{l}. 
        \end{split}
    \end{align*} 
    The above expression is equal to $1$ for the first $26$ gees in \cite[Table 2]{davis2017topological}, since $\phi_{1,2,3}=\psi_1=1$. For the $27$-th gee, it is again $1$, thanks to $\phi_{2,3}\psi_1$.

Now we consider the cases when $k$ is odd, say $k=2l+1$.
Recall the element $\bar{A}_{l}$ that we used in this case of $k=2l$ was $\bar{A}_{l}:= \prod\limits_{j=1}^{l}A_j,$ where $$A_j=\bar{w}_{2j}^{m-1} \bar{u}_{2j}^2\bar{v}_{2j} \bar{R}_{2j}^{m-3}.$$ 
The nonzeroness of the product $\bar{A}_{l}\bar{w}_{2l+1}\bar{u}_{2l+1}\bar{v}_{2l+1}\bar{R}_{2l+1}^{m-3}$ is given as follows.
 Consider the classes $_iY_w:= p^*_i(Y_w)$.
Observe that the type of the $(m,m-1,\dots,m,m-1,m)$-multidegree term of  
$\bar{A}_{l} \bar{w}_{2l+1}\bar{u}_{2l+1}\bar{v}_{2l+1}\bar{R}_{2l+1}^{m-3}$ is $\bar{A}_{l}~ _{2l+1}Y_{1,2,3}$.
Since $\phi_{1,2,3}=1$ in first $26$ cases of \cite[Table 2]{davis2017topological} and $(\phi\otimes \psi)^{\otimes{l}}(\bar{A}_{l})=1$, we get $$((\phi\otimes \psi)^{\otimes{l}}\otimes \phi)(\bar{A}_{l}~ _{2l+1}Y_{1,2,3})=(\phi\otimes \psi)^{\otimes{l}}(\bar{A}_{l})\phi( _{2l+1}Y_{1,2,3})=1.$$
This proves $\bar{A}_{l} \bar{w}_{2l+1}\bar{u}_{2l+1}\bar{v}_{2l+1}\bar{R}_{2l+1}^{m-3}\neq 0$ and thus we get the left inequality of \eqref{eq: tck-type2}. 
For the $27$-th case we consider $\bar{A}_{l} \bar{w}_{2l+1}\bar{u}_{2l+1}^2\bar{R}_{2l+1}^{m-3}$. Clearly, the type of the $(m,m-1,\dots,m,m-1,m)$-multidegree term of $\bar{A}_{l} \bar{w}_{2l+1}\bar{u}_{2l+1}^2\bar{R}_{2l+1}^{m-3}$ is $\bar{A}_{l}~ _{2l+1}Y_{1,2}$. Since $\phi_{1,2}=1$, we obtain that 
$$((\phi\otimes \psi)^{\otimes{l}}\otimes \phi)(\bar{A}_{l}~ _{2l+1}Y_{1,2})=(\phi\otimes \psi)^{\otimes{l}}(\bar{A}_{l})\phi( _{2l+1}Y_{1,2})=1.$$
This gives us the left inequality of \eqref{eq: tck-type2}.

\medskip

\noindent \textbf{Case 2}: \emph{Suppose $m-1$ is a $2$-power and $\alpha$ has doesn't have gees as in last case of \cite[Table 2]{davis2017topological}.}

For $k=2$, Davis proved that $\bar{V}_1^{m} \bar{V}_{2}^2 \bar{V}_3 \bar{R}^{m-4} \neq 0$ by constructing a uniform homomorphism $\psi$, and then applying $\phi \otimes \psi$ to the expansion of the cohomology class in bidegree $(m, m - 1)$, where $\phi$ denotes the Poincar\'e duality isomorphism.
We generalize this idea for general $k\geq 2$.

First, assume that $k$ is even, say $k=2l$.
Consider the element $$\bar{A}_{l}:= \prod\limits_{j=1}^{l}A_j,$$ 
    where $A_j=\bar{w}_{2j}^m \bar{u}_{2j}^2\bar{v}_{2j} \bar{R}_{2j}^{m-4} $.
    We apply $(\phi \otimes \psi)^{\otimes l}$ on the $(m,m-1,m, m-1,\dots, m, m-1)$-multidegree term in the expansion of $\bar{A}_{l}$. 
    As a result, we obtain the following expression: 
$$(\phi \otimes \psi)^{\otimes l} (\bar{A}_{l})= \prod_{j=1}^{l}(\phi\otimes\psi)^{\otimes l}_j (A_j)=(\phi_{1}(\psi_{2,3}+\psi_{1,2,3})+\phi_{1,2,3}\psi_1+\phi_{1,3}\psi_{1,2})^{l}.$$ 
Note that, $\psi_{1,2,3}=\psi_{2,3}=\phi_{1,3}=\phi_{2,3}=0$ and $\phi_{1,2,3}=\psi_1=1$ as explained in the beginning of this section. Hence, the above expression becomes $1$, giving us the desired lower bound of \eqref{eq: tck-type2}.    

Now we consider the cases when $k$ is odd, say $k=2l+1$.
Recall the element $\bar{A}_{l}$ that we used in case $1$ of $k=2l$ was $\bar{A}_{l}:= \prod\limits_{j=1}^{l}A_j,$
    where $$A_j=\bar{w}_{2j}^m \bar{u}_{2j}^2\bar{v}_{2j} \bar{R}_{2j}^{m-4}. $$ We want to show that the product $\bar{A}_{l} \bar{w}_{2l+1}\bar{u}_{2l+1}\bar{v}_{2l+1}\bar{R}_{2l+1}^{m-3}$ is nonzero. Consider the classes $_iY_w:= p^*_i(Y_w)$.
Observe that the type of the $(m,m-1,\dots,m,m-1,m)$-multidegree term of  
$\bar{A}_{l} \bar{w}_{2l+1}\bar{u}_{2l+1}\bar{v}_{2l+1}\bar{R}_{2l+1}^{m-3}$ is $\bar{A}_{l}~ _{2l+1}Y_{1,2,3}$.
Since $\phi_{1,2,3}=1$ and $(\phi\otimes \psi)^{\otimes{l}}(\bar{A}_{l})=1$, we get $$((\phi\otimes \psi)^{\otimes{l}}\otimes \phi)(\bar{A}_{l}~ _{2l+1}Y_{1,2,3})=(\phi\otimes \psi)^{\otimes{l}}(\bar{A}_{l})\phi( _{2l+1}Y_{1,2,3})=1.$$
This proves $\bar{A}_{l} \bar{w}_{2l+1}\bar{u}_{2l+1}\bar{v}_{2l+1}\bar{R}_{2l+1}^{m-3}\neq 0$ and we get the left inequality of \eqref{eq: tck-type2}.

\end{proof}

\section{Monogenic codes of size 4}\label{sec:size4}
   
In this section, we obtain bounds on the higher topological complexity of $\mbalpha$, where the genetic code $\alpha$ is $\langle\{a, a+b, a+b+c,n\}\rangle$ with $n>a+b+c>a+b>a>0$.

Recall that, for a Poincar\'e duality isomorphism $\phi:H^m(\mbalpha;\Z_2)\to \Z_2$, we have
$\phi_0=\phi(R^m)$. 
In \cite[Theorem 4.1]{davis2017topological}, Davis has obtained the expression for $\phi_0$ as follows $$\phi_0=\binom{a}{2}(a+b+c-1)+(a-1)\left(\binom{b}{2}+(b-1)(c-1)\right).$$
The following result classifies values of $a,b,c$ in the genetic code $\langle\{a, a+b,a+b+c,n\}\rangle$ for which the values of $\phi_0$ are odd.

\begin{lemma}\label{lem:phi0-odd}
$\phi_0$ is odd if and only if any of the following conditions hold:
    \begin{enumerate}
        \item $b+c$ even, $a$ even, $a+b~\equiv0 ~( \text{mod}~4)$;
        \item $b+c$ even, $a$ even, $c ~\equiv1 ~( \text{mod}~4)$;
        \item $b+c$ odd and $a~\equiv3 ~(\text{mod}~4)$;
        \item $b+c$ odd, $b~\equiv2,3 ~(\text{mod}~4)$,  and $a~\equiv 0,2 ~(\text{mod}~4)$;
        \item $a,b~\equiv2 ~(\text{mod}~4)$.
    \end{enumerate}    
\end{lemma}

\begin{proof}
    \textit{(1 and 2)}
    Under the assumption that both $a$ and $b+c$ are even, one can write $\phi_0$ as $$ \binom{a}{2}+\binom{b}{2}+ bc + 1 + a(b+c).$$
    By \cite[Lemma 4.2]{davis2017topological}, the last expression is odd under the above-mentioned conditions.

    \textit{(3)} For $a~\equiv3 ~(\text{mod}~4)$,  $(a-1)\left(\binom{b}{2}+(b-1)(c-1)\right) =0$. Moreover, we have $\binom{a}{2}=1$ and $a+b+c-1$ odd. Hence, the result follows. 

    \textit{(4)} For $b+c$ odd, $\phi_0$ becomes $$a \binom{a}{2}+ (a-1)\binom{b}{2}.$$
    When $a~\equiv 0,2 ~(\text{mod}~4)$, it becomes $\binom{b}{2}$, which is odd if and only if $b~\equiv2,3 ~(\text{mod}~4)$.

    \textit{(5)} For $a,b~\equiv2 ~(\text{mod}~4)$, both $\binom{a}{2}$ and $a-1$ are odd. Again, $\binom{b}{2}$ and $b-1$ both are odd. Consequently, there are two odd multiples of $c-1$ and the only surviving odd term is $(a-1)\binom{b}{2}.$
\end{proof}

The previous lemma helps us to obtain the sharp bounds on $\TC_k(\mbalpha)$. 
\begin{theorem}\label{thm:size4-sharpbds}
Let $m$ be a $2$-power and $\langle\{n, a+b+c,a+b, a\}\rangle$ be the genetic code of $\alpha$ satisfying the conditions in \Cref{lem:phi0-odd}. Then
\[km\leq \TC_k(\malpha)\leq km+1.\]
\end{theorem}
\begin{proof}
Note that if $a,b,c$ satisfies conditions in \Cref{lem:phi0-odd}, then $R^m\neq 0$. Then the proof of our assertion is similar to that of \Cref{R^m non-zero}.    
\end{proof}

We denote the intervals $[1,a]$, $(a,a+b]$ and $(a+b, a+b+c]$ by $I_1$, $I_2$ and $I_3$, respectively and classify all possible subgees by tuples of sizes $0$, $1$, $2$ and $3$ (see the proof of \cite[Theorem 4.1]{davis2017topological} for more details).
The tuple $(p)$ of size one represents a subgee containing one element, and the location of this element is decided by the value of $p$. The tuple 
$(p,q)$ of size two represents a subgee $\{i,j\}$ and the locations of $i$ and $j$ are decided by the values of $p$ and $q$, respectively.   
For example, $(1,2)$ refers to a subgee of cardinality two, $\{i,j\}$ such that $i\in I_1$ and $j\in I_2$.
Similarly, the tuple $(p,q,r)$ represents a subgee $\{i,j,k\}$ such that the locations of $i,j,k$ are decided by the values of $p,q,r$. For example, $(1,2,2)$ represents a subgee $\{i,j,k\}$ such that $i\in I_1$ and $j,k\in I_2$.
The following classification of types of elements that can be subgees of the genetic code $\{a,a+b,a+b+c,n\}$ has been given by Davis in the proof of \cite[Theorem 4.1]{davis2017topological}:  
\[\mathcal{S}=\{\emptyset, (1),(2),(3), (1, 1), (1, 2), (1, 3), (2, 2), (2, 3), (1, 1, 1), (1, 1, 2), (1, 1, 3), (1, 2, 2), (1, 2, 3) \}.\]

Next, we aim to obtain the bounds on $\TC_k(\mbalpha)$ when $m>4$. For that purpose, we need subgees of size greater equal $2$.
Such sub-gees are described as follows:
\[
\mathcal{S'}=\{(1, 1), (1, 2), (1, 3), (2, 2), (2, 3), (1, 1, 1), (1, 1, 2), (1, 1, 3), (1, 2, 2), (1, 2, 3)\}.
\]
Now for $U\in \mathcal{S}$, let $u_i$ be the number of $i$'s in $U$ for $1\leq i\leq 3$. For example, if $U=(1,2,2)$, then $u_1=1$ and $u_2=2$.
For $U'\in \mathcal{S'}$, applying $\psi:H^{m-1}(\mbalpha;\Z_2)\to \Z_2$ on the relation $(3)$ of \Cref{thm: cring} we obtain the following expression
\begin{equation}\label{eq: psi-acting-R3}
\sum_{U\in \mathcal{S'}}\binom{a-u_1'}{u_1}\binom{b-u_2'}{u_2}\binom{c-u_3'}{u_3}\psi_U,    
\end{equation}
where $\psi_U=\psi(Y_U)$.

We are now in a position to state our general result.
\begin{theorem}\label{thm:size4general}
With the genetic code of $\alpha$ being as stated above, we have 
    \begin{equation}\label{eq: tck-sec7} 
        km-\bigg\lfloor \frac{k}{2}\bigg\rfloor +1 \leq \TC_k(\mbalpha)\leq km+1
    \end{equation}in the following situations:
    \begin{enumerate}
        \item $m>4$, $a,b \equiv ~1~(\text{mod}~4)$ and $c$ odd;
        \item $m>3$, $a \equiv ~2~(\text{mod}~4)$, $ b  \equiv ~4~(\text{mod}~4)$ and $c$ odd.
    \end{enumerate}
\end{theorem}
\begin{proof}
In the $1$-st situation, following the proof of \cite[Proposition 4.3]{davis2017topological}, we know that there exist a uniform homomorphism $\psi:H^{m-1}(\mbalpha;\Z_2)\to \Z_2$ which sends $Y_{i,j}$ to $1$ and other monomials to $0$. 
Now to achieve our desired assertion,  we will consider two cases depending on whether $m-2$ is not a $2$-power.

Recall that we have $w_j=p_j^*(V_1)$ and $x_j=p_j^*(V_{a+b})$. We now define $y_j:=p_j^*(V_{a+b+c})$.
\noindent\textbf{Case 1}: \emph{Suppose $m-2$ is not a $2$-power.}\\
 Consider the element $$\bar{A}_{l}:= \prod\limits_{j=1}^{l}A_j,$$ 
    where $A_j=\bar{w}_{2j}^{m-2} \bar{x}_{2j}^2\bar{y}_{2j}^{3} \bar{R}_{2j}^{m-4}$.
We deal with the case when $k$ is even, say $k=2l$.  We apply $(\phi \otimes \psi)^{\otimes l}$ on the $(m,m-1,m, m-1,\dots, m, m-1)$-multidegree term in the expansion of $\bar{A}_{l}$. Let us inspect that particular multidegree term of $A_j$, which can get mapped non-trivially. 
    \begin{align*}
        \begin{split}
           &\sum_{i=1}^{m-3} \binom{m-2}{i}\binom{m-4}{m-i-4} {w}_{2j-1}^{i} {x}_{2j-1}^2 {y}_{2j-1}^2 {R}_{2j-1}^{m-i-4} w_{2j}^{m-2-i}y_{2j}R_{2j}^{i}\\
           &+ \sum_{i=1}^{m-3} \binom{m-2}{i}\binom{m-4}{m-i-3} {w}_{2j-1}^{i} {x}_{2j-1}^2 {y}_{2j-1} {R}_{2j-1}^{m-i-3} w_{2j}^{m-2-i}y_{2j}^2R_{2j}^{i-1}  
        \end{split}
    \end{align*}
    All terms of the above sum are of the type $Y_{1,2,3} \otimes Y_{1,3}$. 
    In the first
 situation, it follows from \cite[Theorem 4.1]{davis2017topological} that $\phi$ sends each  $Y_{i,j,k}$ to $1$ and all other monomials to $0$.
 Now applying $\phi \otimes \psi$ to $A_j$, we obtain: 
    \begin{align*}
        \begin{split}
           &\sum_{i=1}^{m-3} \binom{m-2}{i}\binom{m-4}{m-i-4} \phi(Y_{1,2,3}) \psi(Y_{1,3})+ \sum_{i=1}^{m-3} \binom{m-2}{i}\binom{m-4}{m-i-3} \phi(Y_{1,2,3}) \psi(Y_{1,3}) \\
           &= \left[ \sum_{i=1}^{m-3} \binom{m-2}{i}\binom{m-4}{m-i-4}+ \sum_{i=1}^{m-3} \binom{m-2}{i}\binom{m-4}{m-i-3} \right] \phi_{1,2,3}\psi_{1,3}\\
           &=  \binom{2m-6}{m-4} -\binom{m-2}{0}\binom{m-4}{m-4} + \binom{2m-6}{m-3}  = \binom{2m-6}{m-4}+ \binom{2m-6}{m-3}+1  \\
        \end{split}
    \end{align*}
    which is 1, by Lucas's theorem. Therefore, we obtain:
    $$(\phi \otimes \psi)^{\otimes l} (\bar{A}_{l})= \prod_{j=1}^{l}(\phi\otimes\psi)_j  (A_j)= 1.$$

\noindent\textbf{Case 2}: \emph{Suppose $m-2$ is a $2$-power.}\\
    Consider the element $$\bar{A}_{l}:= \prod\limits_{j=1}^{l}A_j,$$ 
    where $A_j=\bar{w}_{2j}^{m-1} \bar{x}_{2j}^2\bar{y}_{2j}^3 \bar{R}_{2j}^{m-5}$.
    Again, we start with $k$ even, say $k=2l$. We inspect the suitable multidegree term of $A_j$, which can get mapped non-trivially under $\phi \otimes \psi$.  
    \begin{align*}
        \begin{split}
           &\sum_{i=1}^{m-2} \binom{m-1}{i}\binom{m-5}{m-i-4} {w}_{2j-1}^{i} {x}_{2j-1}^2 {y}_{2j-1}^2 {R}_{2j-1}^{m-i-4} w_{2j}^{m-1-i}y_{2j}R_{2j}^{i-1}\\
           &+ \sum_{i=1}^{m-2} \binom{m-1}{i}\binom{m-5}{m-i-3} {w}_{2j-1}^{i} {x}_{2j-1}^2 {y}_{2j-1} {R}_{2j-1}^{m-i-3} w_{2j}^{m-1-i}y_{2j}^2R_{2j}^{i-2}  .
        \end{split}
    \end{align*}
    All the terms in the above sum are of the type $Y_{1,2,3} \otimes Y_{1,3}$. Hence, after applying $\phi\otimes \psi$, it becomes, just similar to the last case:
    \begin{align*}
        \begin{split}
           &\sum_{i=1}^{m-2} \binom{m-1}{i}\binom{m-5}{m-i-4} \phi(Y_{1,2,3}) \psi(Y_{1,3})+ \sum_{i=1}^{m-2} \binom{m-1}{i}\binom{m-5}{m-i-3} \phi(Y_{1,2,3}) \psi(Y_{1,3}) \\
           &= \left[ \sum_{i=1}^{m-2} \binom{m-1}{i}\binom{m-5}{m-i-4}+ \sum_{i=1}^{m-2} \binom{m-1}{i}\binom{m-5}{m-i-3}\right] \phi_{1,2,3}\psi_{1,3}\\
           &=  \binom{2m-6}{m-4}   + \binom{2m-6}{m-3}  =   \binom{2m-6}{m-4}  \\
        \end{split}
    \end{align*}
    which is $1$, again by Lucas's Theorem. Therefore, we obtain:
    $$(\phi \otimes \psi)^{\otimes l} (\bar{A}_{l})= \prod_{j=1}^{l}(\phi\otimes\psi)_j (A_j)=  1.
    $$

 For both case-1 and case-2, when $k$ is odd, say $k=2l+1$, we consider the following element $$ \bar{A}_{l} \bar{w}_{2l+1}\bar{x}_{2l+1}^2\bar{y}_{2l+1} \bar{R}_{2l+1}^{m-4}.$$
Then applying $(\phi \otimes \psi)^{\otimes l}\otimes \phi $ on an appropriate multidegree term of the above product, we get $$(\phi \otimes \psi)^{\otimes l}(\bar{A}_{l}) \phi (w_{2l+1}x_{2\ell+1}^2y_{2l+1} R_{2l+1}^{m-4}) = \phi_{1,2,3}=1, $$ and thus completing the proof of the first situation.

Again in the second situation, following the proof of \cite[Proposition 4.4]{davis2017topological} there exists a uniform homomorphism $\psi:H^{m-1}(\mbalpha;\Z_2) \rightarrow \Z_2$ which sends $Y_{i,j}$ to $1$ and other monomials to $0$.
Similarly, as in the first situation, we consider two cases depending on whether $m-2$ is a $2$-power. 

\noindent \textbf{Case 1}: \emph{Suppose $m-2$ is a $2$-power.}\\
    Consider the element $$\bar{A}_{l}:= \prod\limits_{j=1}^{l}A_j,$$ 
    where $A_j=\bar{w}_{2j}^{2} \bar{x}_{2j}^2\bar{y}_{2j}^{m-1} \bar{R}_{2j}^{m-4}$.
    Start with $k$ even, say $k=2l$. We expand $A_j$ in the suitable multidegree, which can get mapped non-trivially under $(\phi \otimes \psi)_j$.
    $$A_j=  \sum_{i=1}^{m-2} \binom{m-1}{i}\binom{m-4}{m-i-2} {x}_{2j-1}^2 {y}_{2j-1}^i {R}_{2j-1}^{m-i-2} w_{2j}^{2}y_{2j}^{m-1-i}R_{2j}^{i-2}.$$
    Note that the above term is of the type $Y_{2,3}\otimes Y_{1,3}$. Thus applying $\phi \otimes \psi$, we get:
    $$\sum_{i=1}^{m-2} \binom{m-1}{i}\binom{m-4}{m-i-2} \phi_{2,3}\psi_{1,3}= \sum_{i=0}^{m-2} \binom{m-1}{i}\binom{m-4}{m-i-2} =\binom{2m-5}{m-2}.$$
    Note that under the $2$-nd situation  we have $\phi_{2,3}=1$ by \cite[Theorem 4.1]{davis2017topological}. As $m-2$ is a $2$-power, the binomial coefficient is of the form 
    $$\binom{2^t +2^{t-1}+ \cdots +2+1}{2^t}.$$
    Hence, by Lucas's theorem, it follows that $\binom{2m-5}{m-2}$ is odd. Therefore, we obtain:
    $$(\phi \otimes \psi)^{\otimes l} (\bar{A}_{l})= \prod_{j=1}^{l}(\phi\otimes\psi)_j (A_j)=  1.$$

\noindent    \textbf{Case 2}: \emph{Suppose $m-2$ is not a $2$-power.}\\
Consider the element $\bar{A}_{l}:= \prod\limits_{j=1}^{l}A_j,$ 
    where $A_j=\bar{w}_{2j}^{2} \bar{x}_{2j}^2\bar{y}_{2j}^{m-2} \bar{R}_{2j}^{m-3}$.

    Start with $k$ even, say $k=2l$. We inspect the suitable multidegree term of $A_j$, which can get mapped non-trivially under $\phi \otimes \psi$. 
    $$A_j=  \sum_{i=1}^{m-3} \binom{m-2}{i}\binom{m-3}{m-i-2} {x}_{2j-1}^2 {y}_{2j-1}^i {R}_{2j-1}^{m-i-2} w_{2j}^{2}y_{2j}^{m-2-i}R_{2j}^{i-1}.$$
    Applying $\phi \otimes \psi$ to the above $Y_{2,3} \otimes Y_{1,3}$ type term we get:
    \begin{align*}
        \begin{split}
            \sum_{i=1}^{m-3} \binom{m-2}{i}\binom{m-3}{m-i-2} \phi_{2,3}\psi_{1,3}&= \sum_{i=0}^{m-2} \binom{m-2}{i}\binom{m-3}{m-i-2}- \binom{m-2}{m-2}\binom{m-3}{0}\\
            &=\binom{2m-5}{m-2}+ 1.
        \end{split}
    \end{align*}
    Again, by Lucas's Theorem, $\binom{2m-5}{m-2}$ is even. Therefore, we obtain:
    $$(\phi \otimes \psi)^{\otimes l} (\bar{A}_{l})= \prod_{j=1}^{l}(\phi\otimes\psi)_j (A_j)=  1.$$   
Now for both case-1 and case-2
 when $k$ is odd, say $k=2l+1$, we consider the following element : $$ \bar{A}_{l}  \bar{x}_{2l+1}^2\bar{y}_{2l+1} \bar{R}_{2l+1}^{m-3}$$
    Then applying $(\phi \otimes \psi)^{\otimes l}\otimes \phi $ on an appropriate multidegree term of the above product, we get: $$(\phi \otimes \psi)^{\otimes l}(\bar{A}_{l}) \phi (x_{2l+1}^2y_{2l+1} R_{2l+1}^{m-3}) = \phi_{2,3}=1. \qedhere$$
\end{proof}

\section{Genetic codes having two genes each of size 3 or genes of Type 1}\label{size3-type1}
In this section, we obtain sharp bounds on the $\TC_k(\mbalpha)$ when the genetic code of $\alpha$ is either having two genes each of size $3$
or having genes of Type $1$.

\subsection{Two genes each of size 3}
In this subsection we inspect the higher topological complexity of $\mbalpha$ where the genetic code of $\alpha$ is $ \langle \{a+b, a+b+c, n\},\{a, a+b+c+d, n\}\rangle$ with $a,b,c,d \geq 1$.

We first show that the $\TC_k(\mbalpha)$ is either $km$ or $km+1$ by classifying values of $a,b,c,d$ of the genetic codes mentioned above for which $\phi_0=1$.  
The expression for $\phi_0$ is given in \cite[Proposition 5.1]{davis2017topological}, which we describe now.
$$ \phi_0= \binom{a-1}{2}+\binom{b}{2}+ bc + (a+1)(b+c+d).$$

\begin{lemma}\label{lem: phi01}
 We have $\phi_0=1$ if and only if either of the following holds
 \begin{enumerate}
     \item $a+b \equiv ~1~(\text{mod}~4)$ or $ a+b+2c  \equiv ~2~(\text{mod}~4)$, and $ (a+1)d \equiv ~0~(\text{mod}~4)$

     \item $a+b \not\equiv ~1~(\text{mod}~4)$ and $ a+b+2c  \not\equiv ~2~(\text{mod}~4)$, and $ (a+1)d \not\equiv ~0~(\text{mod}~4)$
 \end{enumerate}
\end{lemma}

\begin{proof}
    We split the expression of $\phi_0$ into two parts as follows: $$\left[\binom{a-1}{2}+\binom{b}{2}+ bc + (a-1)(b+c)\right] + \left[(a+1)d\right].$$ To have $\phi_0=1$, the two separated terms in the above expression must have different parity. By the if and only if condition given in \cite[Lemma 4.2]{davis2017topological}, the result follows. 
\end{proof}

The proof of the following theorem is similar to that of \Cref{R^m non-zero}, and hence we omit it here.
\begin{theorem} \label{thm-sec7}
Let $m$ be a $2$-power and $ \langle \{a+b, a+b+c, n\},\{a, a+b+c+d, n\}\rangle$ and $a,b,c,d \geq 1$ be the genetic code of $\alpha$. Then for the values of $a,b,c,d$ given in \Cref{lem: phi01}, we have
\[
km\leq \TC_k(\mbalpha)\leq km+1.
\]
\end{theorem}

We now obtain a weaker bound on $\TC_k(\mbalpha)$, generalizing Davis's result from section $5$ of \cite{davis2017topological}.
\begin{theorem}\label{thm-sec7'}
Let $\alpha$ be the genetic code as described above and $m=2^t+m'$ with $2 \leq m'\leq 2^t + 1 $ for some positive integer $t$. Then
\begin{equation}\label{eq: tck-twoofeachsize3}
        km-\bigg\lfloor \frac{k}{2}\bigg\rfloor +1 \leq \TC_k(\mbalpha)\leq km+1.
    \end{equation}
\end{theorem}
\begin{proof} Using similar notations as in \Cref{sec:size4}, we consider the element $$\bar{A}_{l}:= \prod\limits_{j=1}^{l}A_j,$$ 
    where $A_j=\bar{w}_{2j}^{2m'-3} \bar{x}_{2j}^2\bar{y}_{2j} \bar{R}_{2j}^{2^{t+1}-1}$. We first deal with the case when $k$ is even, say $k=2l$.
    Using \cite[Proposition 4.5, Lemma 5.2]{davis2017topological} it follows that $(\phi \otimes \psi)_j(A_j)=1$. Consequently, $\bar{A}_{l}$ is non-zero. For the case where $k=2l+1$, we will consider the element $\bar{A}_{l}\bar{w}_{2l+1} \bar{x}_{2l+1} \bar{R}_{2l+1}^{m-2}$, which will be non-zero. Hence, the inequality \eqref{eq: tck-twoofeachsize3} follows.
\end{proof}

\subsection{Type 1}
In this subsection we inspect the higher topological complexity of $\mbalpha$ where the genetic code $\alpha$ is $ \langle \{1, 1+b, 1+b+c, n\},\{1, 1+b+c+d, n\}\rangle$ with $b,c,d \geq 1$ .

\begin{theorem}\label{thm-sec''}
Let $\alpha$ be the genetic code as described above. Then
\begin{equation}\label{eq: tck-type1}
        km-\bigg\lfloor \frac{k}{2}\bigg\rfloor +1 \leq \TC_k(\mbalpha)\leq km+1
    \end{equation}
    except in either of the following cases:
    \begin{enumerate}
        \item $b \equiv ~1~(\text{mod}~4)$, $c$ odd, and $d$ even
        \item $m-1$ is a $2$-power, or
        \item $m$ is a $2$-power, $b \equiv ~2~(\text{mod}~4)$, $c$ odd and $d$ even.
    \end{enumerate}
\end{theorem}

\begin{proof}
We first define $x_j:=p_j^*(V_{1+b})$ $y_j:=p_j^*(V_{1+b+c})$ following the notations given in \cite[Section 6]{davis2017topological}.
    Now consider the element $\bar{A}_{l}:= \prod\limits_{j=1}^{l}A_j,$ 
    where $A_j=\bar{w}_{2j}^{m-1} \bar{x}_{2j}^2\bar{y}_{2j} \bar{R}_{2j}^{m-3}$. We first deal with the case when $k$ is even, say $k=2l$.
    Using \cite[Proposition 6.2]{davis2017topological} it follows that $(\phi \otimes \psi)_j(A_j)=1$. Consequently, $\bar{A}_{l}$ is non-zero. For the case when $k=2l+1$, we will consider the element $\bar{A}_{l}\bar{w}_{2l+1} \bar{x}_{2l+1}\bar{y}_{2l+1} \bar{R}_{2l+1}^{m-3}$, which will be non-zero. Hence, we obtain the desired inequality \eqref{eq: tck-type1}.
\end{proof}


\noindent\textbf{Acknowledgment.} We thank the anonymous referee for the valuable suggestions which improved the paper in several aspects.
The second author acknowledges the support of NBHM through grant 0204/10/(16)/2023/R\&D-II/2789. The third author acknowledges the support of IISER Pune for the Institute Post-Doctoral fellowship IISER-P/Jng./20235445.

\bibliographystyle{plain} 
\bibliography{references}

\end{document}